\documentclass[12pt]{article}
\newcommand{\R}{\mbox{I$\!$R}}

\newcommand{\qed}{{\hfill {$\rlap{$\sqcap$}\sqcup$}}\\[0.2in]\hspace*{0.5in}}
\newcommand{\qedwh}{{\hfill {$\rlap{$\sqcap$}\sqcup$}}\\[0.2in]}
\newcommand{\bk}{\\[0.03in] \hspace*{0.5in} }
\newcommand{\btd}{\bigtriangledown}

\newcommand{\mfor}{\ \ \ \ {\mbox{for}} \ \ }

\begin{document}

\begin{center} {\LARGE   {\bf Combining Solutions of Semilinear}}
\medskip \smallskip \\ {\LARGE {\bf Partial Differential Equations in ${\R}^n$}} 
\medskip \smallskip \\ {\LARGE {\bf with Critical Exponent}} 
\medskip \medskip  \medskip \\ 
{\large {Man Chun {\Large L}EUNG}}
\medskip \\ {\large National University of Singapore} 
\smallskip \\ {\tt{matlmc@math.nus.edu.sg}}
\end{center}
\vspace{0.1in}
\begin{abstract}
\noindent Let $u_1$ and $u_2$ be two different positive smooth solutions of the equation
$$
\Delta u + n (n - 2) \,u^{{n + 2}\over {n - 2}} = 0 \ \ \ \ {\mbox{in}} \ \ {\R}^n \ \ \ \ (n \ge 3)\,.
$$  
By a result of Gidas, Ni and Nirenberg, 
$$
u_1 (x) = \left( {{\lambda_1}\over {\lambda_1^2 + |x - \xi_1|^2}} \right)^{{n - 2}\over 2}
\ \ \ \ {\mbox{and}} \ \ \ \ 
u_2 (x) = \left( {{\lambda_2}\over {\lambda_2^2 + |x - \xi_2|^2}} \right)^{{n - 2}\over 2}  
$$
for $x \in {\R}^n.$ Here $\xi_1\,, \ \xi_2 \in {\R}^n$ are fixed points and $\lambda_1$ and $\lambda_2$ are positive numbers. Let $u$ be a positive $C^2$-function on ${\R}^n$ such that  $u = u_1$ in $\Omega_1$ and $u = u_2$ in $\Omega_2$, where $\Omega_1$ and $\Omega_2$ are disjoint non-empty open domains in ${\R}^n$. $u$ satisfies the equation 
$$
\Delta u + n (n - 2) K u^{{n + 2}\over {n - 2}} = 0 \ \ \ \ {\mbox{in}} \ \ {\R}^n,
$$  
where 
$$
K (x) := {{- \,\Delta u (x)}\over {n (n - 2) \,u^{{n + 2}\over {n - 2}} (x)}} \mfor x \in {\R}^n.
$$
By the same result of Gidas, Ni and Nirenberg, $K \not\equiv 1$ in ${\R}^n$. In this paper we discuss  lower bounds on $\displaystyle{\sup_{\R^n} |K - 1|}\,.$ Relation with decay estimates  at the isolated singularity via the Kelvin transform is also expounded.
\end{abstract}

\vspace*{0.2in}

2000 AMS MS CLASSIFICATIONS: Primary 35J60\,; Secondary 53C21.

\pagebreak

{\bf \Large {\bf 1. \ \ Introduction}}

\vspace{0.2in}

In this article we consider the quantitative behavior of the scalar curvature equation 
$$
\Delta u + n (n - 2) K u^{{n + 2}\over {n - 2}} = 0 \ \ \ \ {\mbox{in}} \ \ {\R}^n, \ \ n \ge 3\,. \leqno (1.1)
$$
Here $K$ is a smooth function. Equation (1.1) receives extensive studies in connection with the analysis of nonlinear differential operators in ${\R}^n$  and conformal deformations of Riemannian metrics (see, for examples, \cite{Ambrosetti_et_al}, \cite{Chen-Lin.1},  \cite{ChenWX-Li-Ann}, \cite{Kato}, \cite{K-M-P-S}, \cite{Leung.1}, \cite{Leung.5}, \cite{Leung.6}, \cite{Leung.7}, \cite{Lin.1}, \cite{Loewner-Nirenberg}, \cite{Mazzeo-Pacard_2}, \cite{M-P-U}, \cite{Melrose}, \cite{Yanagida_Yotsutani} and the references within). 
We consider the influence of the form of $u$ on $K$. 
Our study is motivated by questions on slow decay for positive solutions of equation (1.1) in ${\R}^n$, when $K$ close to, but not identically equal to, one. In case a solution $u$ fails to have slow decay, it develops bubbles which are modeled after the spherical solutions (see \S 2). We are interested in understanding particular values of $K$ when $u$ is known to be similar to the spherical solutions at certain regions. The investigation leads to estimates on how large $|K(x) - 1|$ can be for some $x$.\bk
Let $u_1$ and $u_2$ be two different positive smooth solutions  of the equation
$$
\Delta u + n (n - 2) \,u^{{n + 2}\over {n - 2}} = 0 \ \ \ \ {\mbox{in}} \ \ {\R}^n.
\leqno (1.2)
$$  
By a salient result of Gidas, Ni and Nirenberg (\cite{Gidas-Ni-Nirenberg.1}, \cite{Gidas-Ni-Nirenberg.2}; cf. \cite{Caffarelli-Gidas-Spruck}),  
$$
u_1 (x) = \left( {{\lambda_1}\over {\lambda_1^2 + |x - \xi_1|^2}} \right)^{{n - 2}\over 2}
\ \ \ \ {\mbox{and}} \ \ \ \ 
u_2 (x) = \left( {{\lambda_2}\over {\lambda_2^2 + |x - \xi_2|^2}} \right)^{{n - 2}\over 2}  
 \leqno (1.3)
$$
for $x \in {\R}^n$ (these are known as spherical solutions). Here $\xi_1\,, \ \xi_2 \in {\R}^n$ are fixed points, and $\lambda_1$ and $\lambda_2$ are positive numbers.   
For a positive $C^2$-function $u$ on ${\R}^n$. define the $K$-function (or curvature function) of $u$ to be 
$$
K (x) := {{- \Delta u (x)}\over {n (n - 2) \,u^{{n + 2}\over {n - 2}} (x)}} \mfor x \in {\R}^n. 
\leqno (1.4)
$$
With the function $K$, $u$ satisfies equation (1.1). 
Suppose that $u = u_1$ in $\Omega_1$ and $u = u_2$ in $\Omega_2$, where $\Omega_1$ and $\Omega_2$ are disjoint non-empty open domains in ${\R}^n$. 
By the same result of Gidas, Ni and Nirenberg, $K \not\equiv 1\,.$ Our task is to provide  lower bounds on $\displaystyle{\sup_{\R^n} |K - 1|}\,.$

\pagebreak

\hspace*{0.5in}In the consideration of slow decay for positive smooth solutions of equation (1.1) in ${\R}^n$, the blow-up analysis involves functions of the form (1.3) (cf. \S 2). 
Assume that $a^2 \le K (x)\le b^2$ for $x \in {\R}^n,\,$ where $a$ and $b$ are positive numbers. One may ask whether when $a/b$ is close to one, positive smooth solutions $u$ of (1.1) have the {\it slow decay}
$$
u (x) \le {C \over {|x|^{(n - 2)/2}}} \mfor |x| \gg 1\,. \leqno (1.5)
$$
Here $C$ is a positive constant. To establish (1.5) is an important first step toward a better understanding on the properties of $u$, which include asymptotic symmetry and completeness (cf. \cite{Chen-Lin.1}, \cite{Chen-Lin.4} and \cite{Leung.8}). See also   \cite{Leung.5} and \cite{Leung.6} for relation of slow decay to asymptotic geometry. Information on slow decay is helpful in finding solutions of equation (1.1), which is also related to the Kazdan-Warner problem \cite{Aubin.1} \cite{Kazdan}.\bk 
By studying addition of spherical solutions, Taliaferro  constructs positive solutions {\it without\,} slow decay and considers the ratio $a/b\,$. This leads him to formulate a conjecture pertaining slow decay and $a/b$ (\cite{Taliaferro.1}, cf. also \cite{Lin.1} and \cite{Taliaferro.2}). Separately, C.-S. Lin conjectures in \cite{Lin.1} that if the Kelvin transform of $K$, given by $K (x/|x|^2)$ for $x \not=0$, extends as a H\"older-continuous function across the origin, then every positive smooth solutions of (1.1) has slow decay. 
Recently the author in \cite{Leung.9} constructs counter-examples to the conjecture of Lin when $n > 4$, and to Taliaferro's conjecture. It appears that $(n - 2)/2$ is a critical order and one may ask the question: Suppose that
$$
|K (x) - 1| \le C |x|^{-{{n - 2}\over 2}} \mfor  |x| \gg 1\,,
$$
where $C$ is a positive constant, does every positive $C^2$-solution $u$ of equation (1.1)  
have slow decay (1.5)?  
The questions on slow decay, together with its applications, are key motivations for us to investigate blending spherical solutions.\bk
In \S 2 we discuss the appearance of bubbles, the blow-up analysis that is essential in the study of equation (1.1) when $K > 0\,.$  The glue-in of a spherical solution is described in the \S 6, together with estimates on the change of the $K$-function. In \S 3 and \S 4, we prove the following main results.\\[0.2in]
{\bf Theorem A.} \ \ {\it Consider a positive function  $u \in C^2 ({\R}^n)$ with}
$$
u (x) = u_1 (x) = \left( {{\lambda_1}\over {\lambda_1^2 + |x|^2}} \right)^{{n - 2}\over 2}
 \ \ \ \ for \ \ |x| < \rho
$$
{\it and}
$$
u (x) = u_2 (x) = \left( {{\lambda_2}\over {\lambda_2^2 + |x|^2}} \right)^{{n - 2}\over 2}
 \ \ \ \ for \ \ |x| \ge R\,.
$$
{\it Here $0 < \rho < R$ are positive numbers. Assume that either} 
$$
{{\lambda_1}\over {\lambda_2}} \le {{\rho^2}\over {R^2}} \left( {1\over {1 + 
\lambda_2^2/R^2}}\right)\,, \ \ \ \ or \ \ \  \
{{\lambda_1^2}\over {\lambda_2^2}} \ge {{3 (n + 2)}\over {2 (n - 2)}} \left(1 + {{R^4}\over {\lambda_2^4}} \right)\,. \leqno (1.6)
$$
{\it Then the $K$-function of $u$ satisfies}
$$
\sup_{x \in B_o (R)} |K (x) - 1| \ge {{n + 2}\over {n}}\,. \leqno (1.7)
$$

\vspace*{0.2in}


\hspace*{0.5in}By using the Kelvin transform, there is a dual condition to (1.6) (remark 4.20). Taking into the account of the ratio $\rho : \lambda_1\,,\,$  the first condition in (1.6) can be interpreted as a relation between the values of $u_1$ and $u_2$ on the boundaries (see remark 4.13).\\[0.2in]
{\bf Theorem B.} \ \ {\it Let $D_1$ and $D_2$ be disjoint open domain in ${\R}^n$ with 
$$
B_{\xi_1} (r_1) \subset D_1 \ \ \ \ and \ \ \ \ B_{\xi_2} (a) \subset D_2\,,
$$ 
where $\xi_1\,, \ \xi_2 \in {\R}^n,$ and $r_1$ and $a$ are positive numbers. Consider a positive $C^2$-function $u$ on ${\R}^n$ with}
$$
u (x) = u_1 (x) = \left( {{\lambda_1}\over {\lambda_1^2 + |x - \xi_1|^2}} \right)^{{n - 2}\over 2} \ \ \ \ for \ \ x \in D_1
$$
{\it and}
$$
u (x) = u_2 (x) = \left( {{\lambda_2}\over {\lambda_2^2 + |x - \xi_2|^2}} \right)^{{n - 2}\over 2}
 \ \ \ \ for \ \ x \in D_2\,.
$$
{\it  Assume that both $u$ and the $K$-function of $u$ have removable singularities at the origin under the Kelvin transformation $x \mapsto x/|x|^2\,,$ $x \not= 0$. 
Given a number $\sigma \ge 1$, if  $\,r_1 \ge \lambda_1\,,$ $a \ge  \,\lambda_2\,$  and}  
$$
{{\lambda_2^2}\over {\lambda_1^2}}  \ge  8^n n \, {{|\xi_1 - \xi_2|^4}\over {r_1^4}} \,
(\sigma^2 + 6)\,, \leqno (1.8)
$$ 
{\it then we have}
$$
\sup_{x \in {\R}^n} |K (x) - 1| \ge {{n + 2}\over {2n}} \, \sigma^2\,. \leqno (1.9)
$$

\vspace*{0.2in}


\hspace*{0.5in}In the last section we apply the idea of estimating $K$-functions to study slow decay (1.5) of positive solutions of equation (1.1). We find restriction on the ``depths" of bubbles.\bk 
Throughout this article $n$ is an integer bigger than two, and $B_x (r)$ is the open ball in ${\R}^n$ with center at $x$ and radius $r > 0\,.$ In conformity with the convention,  we denote  by $C, \ C', \ C_o, \ C_1,...$ various positive constants, which may be different from section to section, unless it is stated otherwise.


\vspace*{0.3in}

{\bf \Large {\bf 2. \ \ The Blow-up Situation}}

\vspace*{0.2in}

In order to discuss the asymptotic behavior of equation (1.1) near infinity, we apply the Kelvin transformation $x \mapsto x/|x|^2$ and consider a positive smooth solution $u$ of the equation
$$
\Delta u +  n (n - 2) K u^{{n + 2}\over {n - 2}} = 0 \ \ \ \ {\mbox{in}} \ \ B_o (1) \setminus \{ 0 \}\,, \leqno (2.1)
$$
where $K$ is a smooth function on $B_o (1) \setminus \{ 0 \}\,.$ 
We examine the isolated singularity of $u$ at $0$. Precisely, assume that $u$ does {\it not} satisfy the slow decay 
$$
u (x) \le {{C_o}\over{ |x|^{(n - 2)/2}}} \mfor x \ \ {\mbox{ with \ \ small}}  \ \  |x|\,,\ \  x \not= 0\,. \leqno (2.2)
$$
\hspace*{0.5in}Given a small positive number $\varepsilon$ and $x \in B_o (5/8) \setminus B_o (\varepsilon)$,\, define
$$
d_\varepsilon (x) := \min\, \{ |x| - \varepsilon\,, \ 5/8 - |x| \}  \ \ \ \ {\mbox{and}} \ \ \ \ U_\varepsilon (x) :=  [d_\varepsilon (x)]^{{n - 2}\over 2} \,u (x)\,.  
$$
Let 
$$
M_\varepsilon := \max_{B_o (5/8)\setminus B_o (\varepsilon)} U_\varepsilon\,.  \leqno (2.3)
$$
Suppose that $U_\varepsilon$ achieves the maximal value 
 $M_\varepsilon$ at a point 
$x_o \in B_o (5/8)\setminus B_o (\varepsilon)\,.$ 
We normalize the solution and spread it out by setting
$$
\lambda := {1\over {[u(x_o)]^{{2\over {n -2}}} }} \ \ \ \ {\mbox{and}} \ \ \ \ w_\lambda (x) := \lambda^{{n - 2}\over 2} u (x_o + \lambda x) = {{u (x_o + \lambda x)}\over {u (x_o)}}\,. \leqno (2.4)
$$
{\it Assume} that 
$$
\lim_{x \to 0} K (x) = 1 \ \ \ \ {\mbox{and}} \ \ \ \ |x| \cdot |\btd K (x)| \le C_1 \mfor x \in B_o (5/8) \setminus \{ 0\}\,. \leqno (2.5)
$$
\hspace*{0.5in}As (2.2) fails to hold, we may assume that $U_\varepsilon (x_o)$ is large when $\varepsilon$ is small.  Bubbles start to emerge in the following sense (see \cite{K-M-P-S} and \cite{Leung.7}, cf. \cite{Cheung-Leung.1}). Given any positive constants $R$ and $\displaystyle{\delta \in (0, 1/2]\,,}$ there is a positive constant $C$ depending on $n$, $R$ and $\delta$ only, such that if $M_\varepsilon \ge C$, then 
$$
\|w_\lambda - u_{\mu, \,y_o} \|_{\,C^2 (B_o (R))} < \delta\,. \leqno (2.6)
$$
Here 
$$
u_{\mu, \,y_o} (x) = \left( {{\mu}\over {\mu^2 + |x - y_o|^2}} \right)^{{n - 2}\over 2}
\mfor x \in {\R}^n
$$
is a spherical solution. From (2.3) and (2.4), we note that  
$$
\lambda = {{d_\varepsilon (x_o)}\over {M_\varepsilon^{2\over {n - 2}}}} < {5\over 8}\, {1\over {M_\varepsilon^{2\over {n - 2}}}} \le {5\over 8} \, {1\over {C^{2\over {n - 2}}}}\,. \leqno (2.7)
$$
\hspace*{0.5in}It is rather clear that there are infinite number of bubbles if (2.2) fails to hold. Indeed, let $C'$ be a positive number such that 
$$
C' > (5/8)^{{n - 2}\over 2}\, u (x_o) > C\,.
$$ 
Here $C$ is the same constant as in (2.7). There exist a small positive number $\varepsilon'$ and $x \in B_o (5/8) \setminus B_o (\varepsilon')$ such that 
$U_{\varepsilon'} (x) \ge C'\,.$ 
The function $U_{\varepsilon'}$ reaches the maximum at a point $x' \in B_o (5/8) \setminus B_o (\varepsilon')$. As
$$
d_{\varepsilon'}^{{n - 2}\over 2} (x') \,u (x') =  U_{\varepsilon'} (x') \ge C' > (5/8)^{{n - 2}\over 2} \,u (x_o)
$$
and $d_{\varepsilon'}  (x') < 5/8\,,$ 
we conclude that $x' \not= x_o$. 
A bubble estimate of the form (2.6) holds near $x'$ for $u$. Let $\lambda' := [u(x')]^{-\,2/( n -2)}$. Given the number $R$, if $C'$ is large enough, then we have $B_{x_o} (\lambda R) \cap B_{x'} (\lambda' R) = \emptyset\,.$\bk
Innate to the construction, we find small flexibility on the choice of the center.   We may  shift  the center to a nearby point $x_1$ (cf. \cite{K-M-P-S}) and obtain
$$
\| u_{1, 0} -  w_\lambda \|_{C^2 (B_o (R))} < \delta\,. \leqno (2.8)
$$
Here $w_\lambda (x)$ becomes $\lambda^{(n - 2)/2} u (x_1 + \lambda x)\,.$ The distance between $x_o$ and $x_1$ is bounded by $c_n\, \lambda$, where $c_n$ is a positive constant that depends on $n$ only (cf. \cite{K-M-P-S}). 
Using the change of variables $y = \lambda \,x$ (without loss of generality, we may assume that $\lambda < 1$), \,(2.8) implies that  
$$
\bigg\Vert \left( {\lambda  \over {\lambda^2 + |y|^2 }} \right)^{{n - 2}\over 2} -  u \,(x_1 + y) \bigg\Vert_{C^o (B_o (\lambda R))} < \delta\,\lambda^{{2 - n}\over 2} \leqno (2.9)
$$
and 
$$
\bigg\Vert \left( {\lambda  \over {\lambda^2 + |y|^2 }} \right)^{{n - 2}\over 2} -  u \,(x_1 + y) \bigg\Vert_{C^1 (B_o (\lambda R))} < {{\delta}\over {\lambda}} \, \lambda^{{2 - n}\over 2}  \leqno (2.10)
$$
for $|y| < \lambda R\,.$ In addition, $R$ can be as large as required.\bk
Under condition (2.5), it is not known to us whether positive smooth solutions of equation (2.1) must have slow decay. For further discussion on the blow-up analysis and slow decay, we refer to \cite{Chen-Lin.2}, \cite{Chen-Lin.3}, \cite{YYLi.4} and \cite{Lin.1}. 

\vspace*{0.3in}

{\bf \Large {\bf 3. \ \ Estimates}}

\vspace*{0.2in}

Let $H (x, \xi)$ be the fundamental solution for the Laplacian on ${\R}^n.\,$ It is given by 
$$
H (x, \xi) = {1\over {(2 - n) \,\omega_n}} {1\over {|x - \xi|^{n - 2}}}
$$
for $x, \ \xi \in {\R}^n$ with $x \not= \xi\,,$ where $\omega_n$ is the volume of the unit sphere in ${\R}^n$. We have $H < 0\,.$\\[0.2in] 
{\bf Theorem 3.1.} \ \ {\it Let $\Omega_1$, $\Omega_2$ and $\Omega$ be non-empty domains in ${\R}^n$ with 
$$
{\overline \Omega}_1 \subset \Omega_2 \subset {\overline \Omega}_2 \subset \Omega\,.
$$ 
Assume that the boundary $\partial \Omega_2$ of $\,\Omega_2$ is smooth. Consider two positive smooth solutions $u_1$ and $u_2$   of the equation 
$$
\Delta u + n (n - 2)\,u^{{n + 2}\over {n - 2}} = 0 \ \ \ \ {\mbox{in}} \ \ \Omega\,. \leqno (3.2)
$$
Let $u_c$ be a positive $C^2$-function on $\Omega$ with $u_c = u_1$ in $\Omega_1$ and $u_c = u_2$ in $\Omega \setminus \Omega_2\,.$ The $K$-function of $u$ is given by  
$$
K (x) := - {{\Delta u_c (x)}\over  { n (n - 2) }} \, [u_c (x)]^{-{{n + 2}\over {n - 2}}} \ \ \ \ for \ \ x \in \Omega\,,
$$
so that  
$$
\Delta u_c + n (n - 2) K u_c^{{n + 2}\over {n - 2}} = 0 \ \ \ \ {\mbox{in}} \ \ \Omega\,.
$$
For any $\xi \in \Omega_1$ we have}
\begin{eqnarray*}
(3.3) \ \ \ \ 
& \,& 4\,n \int_{\Omega_2} H (x, \xi)\, \left[ \,K (x) - 1\, \right] \, dx\\ 
&=& u_1^{-{{4}\over {n - 2}}} (\xi) - u_2^{-{{4}\over {n - 2}}} ( \xi)\\  
& \ & \ \ \ \  + \, (n + 2) \int_{\Omega_2} |H (x, \xi)| 
\left( \bigg\vert \btd u_c^{-{{2}\over {n - 2}}} (x) \bigg\vert^2 - \bigg\vert \btd u_2^{-{{2}\over {n - 2}}} (x) \bigg\vert^2 \right) \,dx\,. \ \ \ \ \ \ \ \ \   
\end{eqnarray*}

\vspace*{0.1in}


{\bf Proof.} We have 
$$
\Delta f (u) = f' (u) \,\Delta u + f'' (u) |\btd u|^2,
$$
where $f$ is a smooth function on $\R^+$.  Hence we obtain
\begin{eqnarray*}
(3.4) \ \ \ \ \ \Delta u_c^{-{{4}\over {n - 2}}}  & = & - {{4}\over {n - 2}} \,u_c^{-{{n + 2}\over {n - 2}}} \,\Delta u_c + \left( {{4}\over {n - 2}} \right) \left( {{n + 2}\over {n - 2}} \right) 
\,u_c^{-{{2n}\over {n - 2}}} \, |\btd u_c|^2 \ \ \ \ \ \ \ \ \\
& = & 4n \,K + (n + 2) \,\bigg\vert \btd u_c^{-{{2}\over {n - 2}}} \bigg\vert^2. 
\end{eqnarray*}
Given $\xi \in \Omega_1 \subset \Omega_2$, it follows from the representation formula (\cite{Gilbarg-Trudinger}, \S 2.4, p. 18) that 
\begin{eqnarray*}
(3.5) \   u_c^{-{{4}\over {n - 2}}} \,(\xi) 
& = & \int_{\Omega_2} H (x, \xi)\, \Delta u_c^{-{{4}\over {n - 2}}} (x) \, dx\\
& \ & \   
 + \int_{\partial \Omega_2} \left[ u_c^{-{{4}\over {n - 2}}} (x) \,{{\partial H (x, \xi)}\over {\partial {\bf n}_x }}\,   - H (x, \xi)\, {{\partial}\over {\partial {\bf n}_x}}  \left( u_c^{-{{4}\over {n - 2}}} (x) \right) \right] \, dS_x\\
& = & 4n \int_{\Omega_2} H (x, \xi)\, K (x) \, dx + (n + 2) \int_{\Omega_2} H (x, \xi) \,\bigg\vert \btd u_c^{-{{2}\over {n - 2}}} (x) \bigg\vert^2 \,dx\\
& \ & \  
+ \int_{\partial \Omega_2} \left[ u_c^{-{{4}\over {n - 2}}} (x) \,{{\partial H (x, \xi)}\over {\partial {\bf n}_x }}\,   - H (x, \xi)\, {{\partial}\over {\partial {\bf n}_x}}  \left( u_c^{-{{4}\over {n - 2}}} (x) \right) \right] \, dS_x  \ \ \ \ \ \ \ \ \ \ \ \\
& = & 4n \int_{\Omega_2} H (x, \xi)\, K (x) \, dx - (n + 2) \int_{\Omega_2} |H (x, \xi)| \,\bigg\vert \btd u_c^{-{{2}\over {n - 2}}} (x) \bigg\vert^2 \,dx\\
& \ & \  
+ \int_{\partial \Omega_2} \left[ u_c^{-{{4}\over {n - 2}}} (x) \,{{\partial H (x, \xi)}\over {\partial {\bf n}_x }}\,   - H (x, \xi)\, {{\partial}\over {\partial {\bf n}_x}}  \left( u_c^{-{{4}\over {n - 2}}} (x) \right) \right] \, dS_x\,,
\end{eqnarray*}
where $\bf n$ is the unit outward normal to $\partial \Omega_2\,.$ Likewise, we can draw the conclusion that 

\newpage

\begin{eqnarray*}
(3.6) \ u_2^{-{{4}\over {n - 2}}} \,(\xi) & = & 4n \int_{\Omega_2} \, H (x, \xi)  \, dx - (n + 2) \int_{\Omega_2} |H (x, \xi)| \,
\bigg\vert \btd u_2^{-{{2}\over {n - 2}}} (x) \bigg\vert^2 \,dx \ \ \ \ \ \ \  \ \ \ \ \ \ \ \\
& \ &  \ + \int_{\partial \Omega_2} \left[ u_2^{-{{4}\over {n - 2}}} (x) \,{{\partial H (x, \xi)}\over {\partial {\bf n}_x }}\,   - H (x, \xi) \, {{\partial}\over {\partial {\bf n}_x}}  \left( u_2^{-{{4}\over {n - 2}}} (x) \right) \right] \, dS_x\,.
\end{eqnarray*}
As $u_c = u_2$ on $\Omega \setminus \Omega_2\,,$ it follows from (3.5) and (3.6)
\begin{eqnarray*}
& \,&u_c^{-{{4}\over {n - 2}}} (\xi) - u_2^{-{{4}\over {n - 2}}} ( \xi)\\ 
&=& 4\,n \int_{\Omega_2} H (x, \xi)\, [ K (x) - 1] \, dx \,  
 - \,(n + 2) \int_{\Omega_2} |H (x, \xi)| \,\bigg\vert \btd u_c^{-{{2}\over {n - 2}}} (x) \bigg\vert^2 \,dx\\
& \ & \ \ \ \  \, + \, (n + 2) \int_{\Omega_2} |H (x, \xi)| \,
\bigg\vert \btd u_2^{-{{2}\over {n - 2}}} (x) \bigg\vert^2 \,dx\,.    
\end{eqnarray*}
Since $u_c (\xi) = u_1 (\xi)$ for $\xi \in \Omega_1\,,$  the above equality leads to  
\begin{eqnarray*}
\ \ \ \ 
& \,& 4\,n \int_{\Omega_2} H (x, \xi)\, [ \,K (x) - 1\,] \, dx\\ 
&=& u_1^{-{{4}\over {n - 2}}} (\xi) - u_2^{-{{4}\over {n - 2}}} ( \xi)\\  
& \ & \ \ \ \  + \, (n + 2) \int_{\Omega_2} |H (x, \xi)| 
\left( \bigg\vert \btd u_c^{-{{2}\over {n - 2}}} (x) \bigg\vert^2 - \bigg\vert \btd u_2^{-{{2}\over {n - 2}}} (x) \bigg\vert^2 \right) \,dx\,. \ \ \ \ \ \ \ \ \   
\end{eqnarray*}
This completes the proof of theorem.\qedwh
{\bf Lemma 3.7.} \ \ {\it Let $\Omega$ be a bounded domain in ${\R}^n.$ Assume that $\Omega \subset B_{x_o} (R)\,.$ We have}
$$
\int_\Omega |H (x, \xi)|\, dx \le {{R^2}\over {2 (n - 2)}} \ \ \ \ for \ \ \xi \in \Omega\,. \leqno (3.8) 
$$

\vspace*{0.1in}

{\bf Proof.} \ \ It suffices to show the case when $\Omega = B_o (R)\,.$ Let $\xi \in B_o (R)$ and $W = B_o (R) \cap B_\xi (R)\,.$ Consider 
$$
{\mbox{I}} := \int_{B_o (R)} |H (x, \xi)|\, dx \ \ \ \ {\mbox{and}} \ \ \ \  {\mbox{II}} := \int_{B_\xi (R)} |H (x, \xi)|\, dx = {{R^2}\over {2 (n - 2)}}\,.
$$
We have $|x - \xi| < R$ for $x \in B_\xi (R) \setminus {\overline W}\,$, \,and $|y - \xi| > R$ for $y \in B_o (R) \setminus {\overline W}\,.$ 
From the form of $H$ and the symmetry of $B_\xi (R) \setminus \overline W$ and $B_o (R) \setminus {\overline W}\!,\,$ one sees that 
\begin{eqnarray*}
{\mbox{II}}  & = & \int_{W} |H (x, \xi)|\, dx + \int_{B_\xi (R) \setminus \overline W} |H (x, \xi)|\, dx\\
& \ge &  \int_{W} |H (x, \xi)|\, dx + \int_{B_o (R) \setminus \overline W} |H (x, \xi)|\, dx = {\mbox{I}}\,,
\end{eqnarray*}
and equality holds only if $\xi = 0$. Thus we obtain the result.\qed
From (3.3) we have
\begin{eqnarray*}
& \,& 4\,n \sup_{x \in \Omega_2} |K (x) - 1| \int_{\Omega_2} |H (x, \xi)|\, dx\\ 
&\ge & u_1^{-{{4}\over {n - 2}}} (\xi) - u_2^{-{{4}\over {n - 2}}} ( \xi)\\  
& \ & \ \ \ \  + \, (n + 2) \int_{\Omega_2} |H (x, \xi)| 
\left( \bigg\vert \btd u_c^{-{{2}\over {n - 2}}} (x) \bigg\vert^2 - \bigg\vert \btd u_2^{-{{2}\over {n - 2}}} (x) \bigg\vert^2 \right) \,dx \ \ \ \ \ \ \ \ \   
\end{eqnarray*}
for $\xi \in \Omega_1\,.$ If $\Omega_2 \subset B_{x_o} (R)\,,$ then it follows from (3.8) that   
\begin{eqnarray*}
(3.9) \ & \  & \sup_{x \in \Omega_2} |K (x) - 1|\\ & \ge & {{n - 2}\over {2n}} R^{-2} \left[ u_1^{-{{4}\over {n - 2}}} (\xi) - u_2^{-{{4}\over {n - 2}}} ( \xi) \right.\\
& \ & \ \ \ \ \   \left. \ + \ (n + 2) \int_{\Omega_2} |H (x, \xi)| 
\left( \bigg\vert \btd u_c^{-{{2}\over {n - 2}}} (x) \bigg\vert^2 - \bigg\vert \btd u_2^{-{{2}\over {n - 2}}} (x) \bigg\vert^2 \right) \,dx \right]\,. \ \ \ \ \ \ \ \ \ \ \ \ \ \ \  
\end{eqnarray*}
In \S 7 we seek to generalize the above result to the case where there is an isolated singularity.

\vspace*{0.3in}

{\bf \Large {\bf 4. \ \ The Concentric Case}}

\vspace*{0.2in}

In (3.9), consider the case when $\Omega_1 = B_o (\rho)$\,, $\Omega_2 = B_o (R),\,$ and 
$$
u_1 (x) = \left( {{\lambda_1}\over {\lambda_1^2 + |x|^2}} \right)^{{n - 2}\over 2}
 \ \ \ \ {\mbox{and}} \ \ \ \ 
u_2 (x) = \left( {{\lambda_2}\over {\lambda_2^2 + |x|^2}} \right)^{{n - 2}\over 2}
\mfor x \in {\R}^n. \leqno (4.1)
$$
Here $\rho$ and $R$ are positive numbers with $0 < \rho < R\,.$ 
In this case we have
$$
u_1^{-{4\over {n - 2}}} (0) - u_2^{-{4\over {n - 2}}} (0) = \lambda_1^2 - \lambda_2^2\,, \leqno (4.2)
$$
and
$$
\bigg\vert \btd u_i^{-{{2}\over {n - 2}}} (x) \bigg\vert^2 = {{4 |x|^2}\over {\lambda_i^2}} \mfor x \in {\R}^n \ \ \ {\mbox{and}} \ \ i = 1, \ 2. \leqno (4.3)
$$
It follows from (4.3) that 

\newpage

\begin{eqnarray*}
& \ &  (n + 2) \int_{\Omega_2} |H (x, 0)| 
\left[ \bigg\vert \btd u_c^{-{2\over {n - 2}}} (x) \bigg\vert^2 - \bigg\vert \btd u_2^{-{2\over {n - 2}}} (x) \bigg\vert^2 \right] \,dx\\ 
& \ge &  \,4 \,(n + 2) \int_{B_o (\rho)} |H (x, 0)| \,  \left( {r^2 \over {\lambda_1^2}} \right) \, dx 
-  \, 4 \,(n + 2) \int_{B_o (R)} |H (x, 0)| \, \left( {r^2\over {\lambda_2^2}} \right) \,dx\\
& = & {{n + 2}\over {n - 2}}\, \left[ {\rho^4 \over {\lambda_1^2}} - {R^4 \over {\lambda_2^2}} \right]\,,
\end{eqnarray*}
where $r =|x|\,.$ we have $K (x) = 1$ for $|x| < \rho$ and 
$$
\int_{B_o (R) \setminus \overline{B_o (\rho)} } |H (x, 0)| \, dx = {{R^2 - \rho^2}\over {2 (n - 2)}}\,.
$$
Using (3.9) we obtain 
\begin{eqnarray*}
(4.4) \ \ \ \ \ \ \ \ \ \ & \ & \sup_{x \in B_o (R) \setminus \overline{B_o (\rho)} } |K (x) - 1|\\
& \ &  \ \ge   {{n - 2}\over {2n}} \left\{ \lambda_1^2 - \lambda_2^2 + {{n + 2}\over {n - 2}}\,  \left[ {\rho^4\over {\lambda_1^2}} - {R^4\over {\lambda_2^2}} \right] \right\} {1\over {R^2 - \rho^2}}\,. \ \ \ \ \ \ \ \ \ \ \ \ \ \ \ \ \ \ \ \ \ \ \ \ \ \ \ \ \   
\end{eqnarray*}
In particular, we have the following results.\\[0.2in]
{\bf Proposition 4.5.} \ \ {\it Let $\lambda_1$, $\lambda_2$, $\rho$ and $R$ be as above. If 
$$
{{\lambda_1}\over {\lambda_2}} \le {{\rho^2}\over {R^2}} \left( {1\over {1 +  \lambda_2^2/R^2}}\right)\,, \leqno (4.6)
$$
then} 
$$
\sup_{x \in B_o (R) \setminus \overline{B_o (\rho)}} |K (x) - 1| \ge {{n + 2}\over {n}}\,.
$$

\vspace{0.1in}

{\bf Proof.} \ \  (4.6) is equivalent to
$$
{{\lambda_1}\over {\rho^2}} \le {{ \lambda_2}\over {R^2}} \left[1  +  \left( {{ \lambda_2}\over R} \right)^2 \right]^{-1}\,. \leqno (4.7)
$$
Hence we have 
$$
{{\rho^4}\over {\lambda_1^2}} \ge R^4 \left[ {{ 1 +  \lambda_2^2 R^{-2} }\over {\lambda_2}} \right]^2 =  
{{R^4} \over {\lambda_2^2}} + 2 R^2 + \lambda_2^2\,.
$$
It follows that 
$$
{{n - 2}\over {2n}} \left\{ \lambda_1^2 - \lambda_2^2 + {{n + 2}\over {n - 2}}\,  \left[ {\rho^4\over {\lambda_1^2}} - {R^4\over {\lambda_2^2}} \right] \right\}  {1\over {R^2 - \rho^2}} \ge  {{n + 2}\over {n}}\,.
$$
This completes the proof.\qedwh
{\bf Proposition 4.8.} \ \ {\it Let $\lambda_1$, $\lambda_2$, $\rho$ and $R$ be as above. If  
$$
{{\lambda_1^2}\over {\lambda_2^2}} \ge {{3(n + 2)}\over {2 (n - 2)}} \left(1 + {{R^4}\over {\lambda_2^4}} \right)\,, \leqno (4.9)
$$
then} 
$$
\sup_{x \in B_o (R) \setminus \overline{B_o (\rho)}} |K (x) - 1| \ge {{n + 2}\over {n}}\,. \leqno (4.10)$$

\vspace{0.1in}

{\bf Proof.} \ \ It can be seen that 
$$
\lambda_2^2  + {{R^4}\over {\lambda_2^2}} \ge 2 R^2. \leqno (4.11)
$$
From (4.9) and (4.11) we have
$$ 
\lambda_1^2 \ge {{3(n + 2)}\over {2 (n - 2)}} \left(  \lambda_2^2  + {{R^4}\over {\lambda_2^2}} \right) 
\ge {{n + 2}\over {n - 2}} \left( \lambda_2^2  + {{R^4}\over {\lambda_2^2}} + R^2 \right)\,. \leqno (4.12)
$$
Hence the result follows from (4.4).\qed
Combining propositions 4.5 and 4.8, we obtain theorem A.\\[0.2in]
{\bf Remark 4.13.} \ \ Let us introduce the depth factors
$$
k_1 := {{\rho}\over {\lambda_1}} \ \ \ \ {\mbox{and}} \ \ \ \ k_2 := {{R}\over {\lambda_2}}\,.
\leqno (4.14)
$$
The depth factor is related to the functions $u_1$ and $u_2$ is the following way. Define
$$
u_i (r) := \left( {{\lambda_i}\over { \lambda_i^2 + r^2}} \right)^{{n - 2}\over 2} \mfor r \ge 0 \ \ {\mbox{and}} \ \ i = 1, \ 2\,.
$$
It can be seen that 
$$
\left[ {{u_1 (0)}\over { u_1 (\rho) }} \right]^{2\over {n - 2}} = {1\over {\lambda_1 [u_1 (\rho)]^{2\over {n - 2}}}} = 1 + k_1^2\,. \leqno (4.15)
$$
Similar expression holds for $k_2\,.$ We see that the smaller the value of $u_1 (\rho)$ (or the deeper the bubble goes), the larger the depth factor. \bk 
(4.15) helps to shed light on condition (4.6) in proposition 4.5. We assert that condition (4.6) is equivalent to the condition that   
$$
u_2 (R) \ge \nu \,u_1 (\rho)\,, \leqno (4.16)
$$
where $\nu$ is the number given by 
$$
\nu := \left( {{k_1^2 + 1}\over {k_1^2}}  \right)^{{n - 2}\over 2}. \leqno (4.17)
$$
Observe that $\nu \to 1$  as $k_1  \to \infty\,.$\bk  
Condition (4.6) can be written as
$$
{{\lambda_2}\over {\lambda_1}} \le  {{k_1^2}\over {k_2^2 + 1}}\,. \leqno (4.18)
$$
Condition (4.16) is the same as 
$$
 \left[ {{u_2 (R)}\over { u_1 (\rho) }} \right]^{2\over {n - 2}} {1\over {\nu^{2\over {n - 2}} }} \ge 1\,. \leqno (4.19)
$$
From (4.15) and (4.19), we obtain  
$$
{{k_1^2}\over {k_2^2 + 1}} =  {{k_1^2 + 1}\over {k_2^2 + 1}} {1\over {\nu^{2\over {n - 2}} }} =  {{\lambda_2}\over {\lambda_1}} \left[ {{u_2 (R)}\over { u_1 (\rho) }} \right]^{2\over {n - 2}} {1\over {\nu^{2\over {n - 2}} }} \ge {{\lambda_2}\over {\lambda_1}}\,.
$$
That is, condition (4.6) is satisfied if $u_2 (R) \ge \nu u_1 (\rho)\,,$ and vice versa.  
We note that the term $\nu u_1 (\rho)$ is independent on $u_2\,.$\\[0.2in]
{\bf Remark 4.20.} \ \ By applying the Kelvin transformation $x \mapsto \rho^2 x/|x|^2$ on $u$, we come to a situation where there is a positive $C^2$-function $\tilde u$ such that 
$$
\tilde u (x) = \tilde u_2 (x) := \left( {{\tilde \lambda_2}\over {{\tilde \lambda}_2^2 + |x|^2}} \right)^{{n - 2}\over 2}
\mfor |x| < \rho^2/R
$$
and 
$$
\tilde u (x) = \tilde u_1 (x) = \left( {{\tilde \lambda_1}\over {{\tilde \lambda}_1^2 + |x|^2}} \right)^{{n - 2}\over 2}
\mfor |x| \ge \rho\,.
$$
Here 
$$
\tilde \lambda_i = \rho^2/\lambda_i \ \ \ \ \mfor \ \ i = 1, \ 2\,.
$$ 
[See \cite{Leung.7}. cf. also (5.5), (5.6) and (5.7) in \S  5; here $\tilde u$ has a removable singularity at the origin because of the form of $u_2$\,.] 
From theorem A, assuming that either
$$
{{\lambda_1}\over {\lambda_2}} \le {{\rho^2}\over {R^2}} \left( {1\over {1 +  \rho^2/ \lambda_1^2}}\right) \,, \ \ \ \ {\mbox{or}}  \ \ \ \
{{\lambda_1^2}\over {\lambda_2^2}} \ge {{3(n + 2)}\over {2 (n - 2)}} \left(1 + {{\lambda_1^4}\over {\rho^4}} \right)\,, \leqno (4.21)
$$
we have  
$$
\sup_{x \in B_o (R) \setminus \overline{B_o (\rho)}} |K (x) - 1| \ = 
\sup_{x \in B_o (\rho) \setminus \overline{B_o (\rho^2/R)}} |\tilde K (x) - 1| \ge {{n + 2}\over {n}}\,.
$$
Here $\tilde K$ is the $K$-function of $\tilde u$ and is given by $\displaystyle{ \tilde K (x) = K (\rho x/|x|^2)}\,$ for $x \not= 0\,$. Thus we may replace the conditions (1.6) by (4.21) in order to obtain the estimate (1.7).\\[0.2in]
\hspace*{0.5in}We note that the results can be generalized into the situation when $B_o (\rho) \subset \Omega_1$ and $\Omega_2 \subset B_o (R)\,,$ and similar technique can be applied to the eccentric case. 
 It is clear that, from (4.4), if the parameters $\rho$, $R$, $\lambda_2$ are fixed, then 
$$
\sup_{x \in B_o (R) \setminus \overline{B_o (\rho)}} |K (x) - 1| \to \infty
\ \ \ \ {\mbox{as}} \ \ \ \  
\lambda_1 \to 0^+.
$$

\vspace*{0.3in}

{\bf \Large {\bf 5. \ \ Disjoint Domains}}

\vspace*{0.2in}

By a translation, in theorem B, we may take $\xi_2 = 0\,.$ Thus 
let $D_1$ and $D_2$ be disjoint open domains in ${\R}^n$ with 
$$
B_o (a) \subset D_2 \ \ \ \ {\mbox{and}} \ \ \ \ B_{\xi_1} (r_1) \subset D_1 
$$
for positive number $r_1$ and $a,$ and for a point $\xi_1 \in {\R}^n$. Consider spherical solutions 
$$
u_1 (x) = \left( {{\lambda_1}\over {\lambda_1^2 + |x - \xi_1|^2}} \right)^{{n - 2}\over 2}
\ \ \ \ {\mbox{and}} \ \ \ \ 
u_2 (x) = \left( {{\lambda_2}\over {\lambda_2^2 + |x|^2}} \right)^{{n - 2}\over 2}  
\mfor x \in {\R}^n,
$$
where $\lambda_1$ and $\lambda_2$ are positive numbers. Let $u_c$ be a positive $C^2$-function on ${\R}^n$ such that 
$$
u_c = u_1 \ \ {\mbox{in}} \ \ D_1  \ \ \ \ {\mbox{and}} \ \ \ \ u_c = u_2 \ \ {\mbox{in}} \ \ D_2\,.
$$ 
$u_c$ satisfies the equation 
$$
\Delta u_c + n (n - 2) \,K u_c^{{n + 2}\over {n - 2}} = 0 \ \ \ \ {\mbox{in}} \ \ {\R}^n, \leqno (5.1)
$$
where $K$ is the $K$-function of $u_c$ (cf. (1.4)). We use a Kelvin transform to bring  the setting into that of theorem 3.1. (cf. also example 5.25.)\bk 
The Kelvin transform of $u_c$ about the sphere of radius $a$ and center at the origin is given by 
$$
{\tilde u}_c (x) := \left( {a\over {|x|}} \right)^{n - 2} u_c \left( {{a^2 x}\over {|x|^2}} \right) \mfor x \in {\R}^n \setminus \{ 0 \}\,. \leqno (5.2)
$$
We find that ${\tilde u}_c$ satisfies the equation \cite{Leung.7}
$$
\Delta {\tilde u}_c + n (n - 2) \,\tilde K {\tilde u}_c^{{n + 2}\over {n - 2}} = 0 \ \ \ \ {\mbox{in}} \ \ {\R}^n \setminus \{ 0 \}\,, \leqno (5.3)
$$
where
$$
\tilde K (x) := K \left( {{a^2 x}\over {|x|^2}} \right) \mfor x \in {\R}^n \setminus \{ 0 \}\,.
$$

\vspace*{0.2in}

{\bf Lemma 5.4.} \ \ {\it Let $u$ be a function defined on ${\R}^n$, and let
$$
\tilde u (x) := {1\over {|x|^{n - 2}}}\, u \left( {x\over {|x|^2}} \right) \ \ \ \ for \ \ x \not= 0 
$$
be the Kelvin transforms of $u$ with center at the origin and radius one; 
and 
$$
\tilde u_{\xi_2, a} (x) := {{a^{n - 2}}\over {|x - \xi_2|^{n - 2}}} \,\,u \left(\xi_2 + {{a^2 (x - \xi_2)}\over {|x - \xi_2|^2}} \right) \ \ \ \ for \ \ x \not= \xi_2\,,
$$
be the Kelvin transforms of $u$ with center at $\xi_2$ and radius $a > 0\,.$  Then}
$$
\tilde u_{\xi_2, a} (x) = {{a^{n - 2}}\over {|x - \xi_2|^{n - 2}}}   \,\, \bigg\vert \xi_2 + {{a^2 (x - \xi_2)}\over {|x - \xi_2|^2}} \bigg\vert^{2 - n} \, \tilde u \left( {{\xi_2 + {{a^2 (x - \xi_2)}\over {|x - \xi_2|^2}}}\over {\bigg\vert \xi_2 + {{a^2 (x - \xi_2)}\over {|x - \xi_2|^2}} \bigg\vert^2}} \right)
$$
{\it  for $x \not= \xi_2\,.$}

\vspace*{0.1in}

{\bf Proof.} \ \ We have
$$
u (x) = {1\over {|x|^{n - 2}}}\, \tilde u \left( {x\over {|x|^2}} \right) \mfor x \not= 0\,.
$$
Substitute into the expression for $\tilde u_{\xi_2, a}$ we obtain the result.\qed
It follows from lemma 5.4 that, as $x \to \xi_2$, 
$$
\tilde u_{\xi_2, a} (x)  \to  {{a^{n - 2}}\over {|x - \xi_2|^{n - 2}}}   \, \bigg\vert  {{a^2 (x - \xi_2)}\over {|x - \xi_2|^2}} \bigg\vert^{2 - n} \, \tilde u \left( {{ {{a^2 (x - \xi_2)}\over {|x - \xi_2|^2}}}\over {\bigg\vert  {{a^2 (x - \xi_2)}\over {|x - \xi_2|^2}} \bigg\vert^2}} \right)
= {1\over {a^{n - 2}}} \,\tilde u \left( {{x - \xi_2}\over {a^2}} \right)\,.
$$
By the assumption in theorem B and lemma 5.4, mindful also of the translation,   ${\tilde u}_c$ can be extended as a $C^2$-function across the origin. Likewise, $\tilde K$ can be extended as a continuous function, so that equation (5.3) is satisfied on the whole ${\R}^n.$ With all this, compare also with theorem 7.9.\bk
The Kelvin transforms of $u_1$ and $u_2$ about the sphere of radius $a > 0$ and center at the origin are given by 
$$
{\tilde u}_1 (x) := \left( {{ {{\overline \lambda}_1}}\over {{\overline \lambda}_1^2 + |x - {\overline \xi}_1|^2}} \right)^{{n - 2}\over 2}
\mfor x \in {\R}^n, \leqno (5.5)
$$
and
$$			
{\tilde u}_2 (x) := \left( {{a^2 \lambda_2^{-1}}\over {a^4 \lambda_2^{-2} + |x|^2}} \right)^{{n - 2}\over 2}  
\mfor x \in {\R}^n, \leqno (5.6)
$$
respectively  (see, for example, \cite{Leung.7}). Here
$$
{\overline \xi}_1 = {{a^2 \,\xi_1} \over { \lambda_1^2 + |\xi_1|^2}} \ \ \ \ {\mbox{and}} \ \ \ \ {\overline \lambda}_1 = {{a^2\lambda_1}\over { \lambda_1^2 + |\xi_1|^2}}\,. \leqno (5.7)
$$
Let $\Omega_1$ denote the image of $D_1$ under the inversion 
$$
x \mapsto a^2 x/|x|^2, \ \ \ \ x\not=0\,.
$$   
It follows that 
$$
{\tilde u}_c = {\tilde u}_1 \ \ \ \ {\mbox{in}} \ \ \ \ \Omega_1 \ \ \ \ {\mbox{and}} \ \  {\tilde u}_c = {\tilde u}_2 \ \ \ \ {\mbox{in}} \ \ {\R}^n \setminus B_o (a)\,. 
$$ 
We seek a suitable point to apply (3.9).\bk
The inversion $x \mapsto a^2 x/|x|^2$ sends a ball to a ball. Consider the line passing through the origin and $\xi_1$. It can be seen that the ``inverted" ball of $B_{\xi_1} (r_1)$ has diameter $d$ given by  
$$
d =  \left( {{a^2}\over {|\xi_1| - r_1}} - {{a^2}\over {|\xi_1| + r_1}} \right) = {{2 r_1\, a^2}\over {|\xi_1|^2 - r_1^2}}\,. \leqno (5.8)
$$ 
($|\xi_1| > r_1$ because $B_o (a)$ and $B_{\xi_1} (r_1)$ are disjoint.)  
{\it Assume} that 
$$
r_1  \ge \lambda_1\,.
$$ 
We have 
$$
{{a^2 |\xi_1|}\over {|\xi_1|^2 + r_1 |\xi_1|}} \,\le \,{{a^2 |\xi_1|}\over {|\xi_1|^2 + \lambda_1^2}} \,\le\, {{a^2 |\xi_1|}\over {|\xi_1|^2 - r_1 |\xi_1|}}\,.
$$
Hence there exists a point $\xi \in B_{\xi_1} (r_1)$ such that 
$$
{{a^2 \xi}\over { |\xi|^2 }} = {{a^2 \xi_1}\over {|\xi_1|^2 + \lambda_1^2}} = {\bar \xi}_1\,. \leqno (5.9)
$$
From (5.5) and (5.9) we obtain 
$$
{\tilde u}_1 (a^2 \xi/ |\xi|^2)  = \left( {1\over { 
{\overline \lambda}_1  }} \right)^{{n - 2}\over 2}\,. 
$$
That is, 
$$
{\tilde u}_1^{{-4}\over {n - 2}} (a^2 \xi/ |\xi|^2) = {\overline \lambda}_1^2 = {{a^4 \lambda_1^2}\over {(\lambda_1^2 + |\xi_1|^2)^2}}. \leqno (5.10) 
$$
Moreover,  
$$
{\tilde u}_2^{{-4}\over {n - 2}} (a^2 \xi/ |\xi|^2) = {{(\lambda_2^{-2} |\xi|^2 + 1)^2 a^4 \lambda_2^2}\over {|\xi|^4}}= \left({a\over {|\xi|}}\right)^4 {{(\lambda_2^2 + |\xi|^2)^2}\over { \lambda_2^2}}\,. \leqno (5.11)
$$
The following lemma can be verified by direct calculations.\\[0.2in]
{\bf Lemma 5.12.} \ \ {\it Given a point $\xi \in {\R}^n$ and a positive number $\lambda$, \,let} 
$$
u  (x) = \left( {{\lambda}\over {\lambda^2 + |x - \xi|^2}} \right)^{{n - 2}\over 2} \ \ \ \ for \ \  x \in {\R}^n.
$$
{\it We have}
$$
\bigg\vert \btd u^{-{{2}\over {n - 2}}} (x) \bigg\vert^2 = 4 \,{{|x - \xi|^2}\over {\lambda^2}} \ \ \ \ for \ \  x \in {\R}^n. \leqno (5.13)
$$

\vspace*{0.2in}
\hspace*{0.5in}Thus, 
$$
\bigg\vert \btd {\tilde u}_1^{-{2\over {n - 2}}} (x)\bigg\vert^2  = 
4 \,{{|x - {\overline \xi}_1|^2}\over {{a^4 \lambda}_1^2}} (\lambda_1^2 + |\xi_1|^2)^2 
\ \ \ \ {\mbox{and}} \ \ \ \ \bigg\vert\btd {\tilde u}_2^{-{2\over {n - 2}}} (x)\bigg\vert^2  = 4 \,|x|^2 \lambda_2^2 \,a^{-4} 
$$
for $x \in {\R}^n.\,$ 
In this case $\Omega_2 = B_o (a)\,.$ It follows from lemma 3.7 that 
\begin{eqnarray*}
(5.14) \ \ \ \ \ \ \ \ & \ & (n + 2) \int_{\Omega_2} |H (x, a^2\,\xi/|\xi|^2)|\, \bigg\vert\btd {\tilde u}_2^{-{2\over {n - 2}}} (x) \bigg\vert^2 dx\\ 
& \le & 4 \,(n + 2)\, \lambda_2^2 \, a^{-2}   \int_{B_o (a)} |H (x, a^2\,\xi/|\xi|^2)|\,dx \le \,2  \left({{n + 2}\over {n - 2}}  \right) \lambda_2^2\,.  
\ \ \ \ \ \ \ \ \ \ \ \ \ \ \ \ \ \ \   \ \  
 \ \  
\end{eqnarray*} 
The inverted ball of $B_{\xi_1} (r_1)$ inside $\Omega_1$ contains a ball $B_{\xi_o} (\rho)$ with radius $\rho = d/8\,.$ By choosing the center $\xi_o$ suitably [cf. (5.8) and (5.9)], we have 
$$
\vert x - a^2 \xi/|\xi|^2\vert = |x - {\overline \xi}_1| \ge \rho \ \ \ \ {\mbox{for \ \ all}} \ \  x \in B_{\xi_o} (\rho)\,.\,
$$ 


Hence
\begin{eqnarray*}
(5.15) \ \ \ \ \ \ \ \ \ \ \ \ \ \
& \ & (n + 2) \int_{\Omega_1} |H (x, a^2\,\xi/|\xi|^2)|\, |\btd {\tilde u}_c^{-{2\over {n - 2}}} (x)|^2 dx\\ 
& \ge & (n + 2) \int_{B_{\xi_o} (\rho)} |H (x, a^2\,\xi/|\xi|^2)|\, |\btd {\tilde u}_1^{-{2\over {n - 2}}} (x)|^2 dx\\
& \ge & 4 \,(n + 2)\, {{(\lambda_1^2 + |\xi_1|^2)^2}\over {a^4 \lambda_1^2}}  
\int_{B_{\xi_o} (\rho)} |H (x, a^2\,\xi/|\xi|^2)|\, |x - {\overline \xi}_1|^2  \, dx\ \ \ \ \ \ \ \ \ \ \ \ \ \ \ \ \ \ \ \ \ \ \ \ \\
& \ge &  4 \,(n + 2)\,  {{(\lambda_1^2 + |\xi_1|^2)^2}\over {a^4 \lambda_1^2}}  \rho^2 \int_{B_{\xi_o} (\rho)} |H (x, a^2\,\xi/|\xi|^2)|  \, dx\\
& \ge &  4(n + 2) \,  {{(\lambda_1^2 + |\xi_1|^2)^2}\over {a^4 \lambda_1^2}}  \rho^2  \left( {1\over {(n - 2) \,\omega_n}} {1\over {d^{n - 2}}} \right) \int_{B_{\xi_o} (\rho)} dx\\
& \ge &  4(n + 2) \,  {{(\lambda_1^2 + |\xi_1|^2)^2}\over {a^4 \lambda_1^2}}  \rho^2  \left( {1\over {(n - 2) \,\omega_n}} {1\over {d^{n - 2}}} \right) \left( {\omega_n\over n} \rho^n\right)\\
& \ge & {{4(n + 2)}\over {n (n - 2)}} {1\over {8^{n + 2}}} \,  {{(\lambda_1^2 + |\xi_1|^2)^2}\over {a^4 \lambda_1^2}}  \,d^4\\
& = & {{4(n + 2)}\over {n (n - 2)}} {1\over {8^{n + 2}}} \,  {{(\lambda_1^2 + |\xi_1|^2)^2}\over {a^4 \lambda_1^2}} \left( {{2^4 \,r_1^4\, a^8}\over{ (|\xi_1|^2 - r_1^2)^4 }}\right)\\
& \ge & {{n + 2}\over {n (n - 2)}} {1\over {8^n}} {{a^4}\over {|\xi_1|^8}} {{(\lambda_1^2 + |\xi_1|^2)^2}\over {\lambda_1^2}} \,r_1^4\,.
\end{eqnarray*}
It follows from (3.9), (5.10), (5.11), (5.14) and (5.15) that 
\begin{eqnarray*}
(5.16) \ \ & \ & \sup_{{\R}^n \setminus B_o (a)} |K  - 1| = \sup_{ B_o (a)} |\tilde K  - 1|\\
& \ge &  
{{n - 2}\over {2n}} \, a^{-2} \left[ \, {{a^4 \lambda_1^2}\over {(\lambda_1^2 + |\xi_1|^2)^2}} 
 + 
{{n + 2}\over {n (n - 2)}} {1\over {8^n}}  {{a^4}\over {|\xi_1|^8}} {{(\lambda_1^2 + |\xi_1|^2)^2}\over { \lambda_1^2}}\, r_1^4 \right.\\
& \ & \ \ \ \ \ \ \ \ \ \ \ \ \ \ \ \ \ \ \ \ \left. - \left({a\over {|\xi|}}\right)^4 {{(\lambda_2^2 + |\xi|^2)^2}\over { \lambda_2^2}}
 - 2\, \left( {{n + 2}\over {n - 2}} \right) \lambda_2^2 \right]\,. \ \ \ \ \ \ \ \ \ \ \ \ \ \ \ \ \ \ \ \ \ \ \ \ \ \ 
\end{eqnarray*}
\hspace*{0.5in}Let the positive numbers $c$, $k$ and $C$ be defined by the equations 
$$
r_1 = c\, \lambda_1\,, \ \ \ \ a = k \,\lambda_2\,, \ \ \ \  |\xi_1| = C \lambda_2\,, \leqno (5.17)
$$
respectively. We obtain 
$$
{{a^4 \lambda_1^2}\over {(\lambda_1^2 + |\xi_1|^2)^2}} = {{k^4 \lambda_2^4 \lambda_1^2}\over {(\lambda_1^2 + C^2 \lambda_2^2)^2}} = {{k^4 \lambda_1^2}\over {\left( {{\lambda_1^2}\over {\lambda_2^2}} + C^2 \right)^2 }}\,,
$$
$$
{{a^4}\over {|\xi_1|^8}} {{(\lambda_1^2 + |\xi_1|^2)^2}\over { \lambda_1^2}}\, r_1^4 = 
  {{k^4 \lambda_2^4}\over {C^8 \lambda_2^8}} \,\left(\lambda_1^2 + C^2 \lambda_2^2\right)^2\, c^4 \lambda_1^2 >   {{c^4}\over {C^4}} k^4 \lambda_1^2\,.
$$
From (5.9) we have
$$
|\xi| = |\xi_1| + {{\lambda_1^2}\over {|\xi_1|}} \ \ \Longrightarrow \ \  |\xi| > |\xi_1|\,.
$$
{\it Assume} that $k = a/\lambda_2 \ge 1\,.$ It follows that  
$$
\left({a\over {|\xi|}}\right)^4 {{(\lambda_2^2 + |\xi|^2)^2}\over { \lambda_2^2}}  = \left( {{\lambda_2^2}\over {|\xi|^2}} + 1 \right)^2 k^4 \lambda_2^2 
\le \left( {{\lambda_2^2}\over {|\xi_1|^2}} + 1 \right)^2 k^4 \lambda_2^2 \le 2^2 k^4 \lambda_2^2\,,
$$
as $|\xi_1| > a \ge \lambda_2\,.$ From (5.16) and the fact that $a^2 = k^2 \lambda_2^2$, we have 
\begin{eqnarray*}
(5.18) \ \ \ \ \ \ \ \ \ \ \ \ \   & \ & \sup_{{\R}^n \setminus B_o (a)} |K  - 1| = \sup_{ B_o (a)} |\tilde K  - 1|\\
& \ge &  
{{n - 2}\over {2n}} \left[ {{k^2}\over {\left( {{\lambda_1^2}\over {\lambda_2^2}} + C^2 \right)^2 }} {{\lambda_1^2}\over {\lambda_2^2}}
+ {{n + 2}\over {n (n - 2)}} {1\over {8^n}}  k^2 {{\lambda_1^2}\over {\lambda_2^2}} {{c^4}\over {C^4}} \right.\\
& \ & \ \ \ \ \ \ \ \ \ \ \  \ \ \ \ \left. -\, \, 4\, k^2 \,- \,2 \left( {{n + 2}\over {n - 2}} \right) {1\over {k^2}} \right]\,.\ \ \ \ \ \ \ \ \ \ \ \ \ \ \ \ \ \ \ \ \ \ \ \ \ \ \ \ \ \ \ \ \ \ 
\end{eqnarray*}
Observe that 
$$
{{k^2}\over {\left( {{\lambda_1^2}\over {\lambda_2^2}} + C^2 \right)^2 }} {{\lambda_1^2}\over {\lambda_2^2}} \,\le \,{{k^2}\over {4 \,C^2}} = {{a^2}\over {4 |\xi_1|^2}} \le {1\over 4}\,.
$$
\hspace*{0.5in}Given a number $\sigma \ge 1$, we seek conditions for the inequality
$$
{{n + 2}\over {n (n - 2)}} {1\over {8^n}}  k^2 {{\lambda_1^2}\over {\lambda_2^2}} {{c^4}\over {C^4}} -\, \, 4\, k^2 \,- \,2 \left( {{n + 2}\over {n - 2}} \right) {1\over {k^2}} \ge \left( {{n + 2}\over {n - 2}} \right) \sigma^2  \leqno (5.19)
$$
to hold. That is, 
$$
{{n + 2}\over {n (n - 2)}} {1\over {8^n}}  {{\lambda_1^2}\over {\lambda_2^2}} {{c^4}\over {C^4}} \ge 4 + 
 \left( {{n + 2}\over {n - 2}} \right) \left[{2\over {k^4}} + {{\sigma^2}\over {k^2}} \right]\,. 
\leqno (5.20)
$$
As
$$
{{c^4}\over {C^4}} = {{\lambda_2^4}\over {\lambda_1^4}} {{r_1^4}\over {|\xi_1|^4}}\,,
$$
(5.20) is equivalent to 
$$
{{\lambda_2^2}\over {\lambda_1^2}} \ge 8^n {{|\xi_1|^4}\over {r_1^4}} \left[    
{{4 n(n - 2)}\over {n + 2}}   + 
   \left({{2n \lambda_2^4}\over {a^4}} + {{n \sigma^2 \lambda_2^2}\over {a^2}} \right)  
\right]\,. \leqno (5.21)
$$
As $\lambda_2 /a \le 1$, (5.21) holds if
$$
{{\lambda_2^2}\over {\lambda_1^2}}  \ge 8^n  n \,{{|\xi_1|^4}\over {r_1^4}}  \left[ 4 + 2 + \sigma^2 \right]\,. \leqno (5.22)
$$
Hence we conclude that if 
$$
{{\lambda_2^2}\over {\lambda_1^2}}  \ge  8^n n \,{{|\xi_1|^4}\over {r_1^4}}   (\sigma^2 + 6)\,, \leqno (5.23)
$$
and $r_1 \ge \lambda_1$ and $a \ge \lambda_2$, then 
$$
\sup_{ {\R}^n } |K - 1| \ge {{n + 2}\over {2n}} \sigma^2.  \leqno (5.24)
$$
Taking into the account of the translation, we obtain theorem B.
\vspace*{0.2in}

{\bf Example 5.25.} \ \ Let 
$$
u (x) := u_1 (x) + u_2 (x)  = \left( {{\lambda_1}\over {\lambda_1^2 + |x - \xi_1|^2}} \right)^{{n - 2}\over 2} + \left( {{\lambda_2}\over {\lambda_2^2 + |x|^2}} \right)^{{n - 2}\over 2} \mfor x \in {\R}^n. 
$$
In can be seen that \cite{Leung.6} 
$$
| K (x) - 1|\, \le \,  1 - 2^{4/(2-n)} \mfor x \in {\R}^n. \leqno (5.26)
$$
Here $K$ is the $K$-function of $u\,.$  At points $x \in {\R}^n$ with $u_1 (x) = u_2 (x)$, we have $K (x) = 2^{4/(2-n)}\,.$ Moreover, 
$$
\lim_{|x| \to \infty} K (x) = {{\lambda_1^{{n + 2}\over 2} + \lambda_2^{{n + 2}\over 2}}\over 
{ \left( \lambda_1^{{n - 2}\over 2} + \lambda_2^{{n - 2}\over 2} \right)^{{n + 2}\over {n - 2}} }}\,. \leqno (5.27)
$$
In particular, if $\lambda_1 = \lambda_2$, we have $\displaystyle{ \lim_{|x| \to \infty} K (x) = 
2^{4/(2-n)}\,.}$\bk
We observe the decay behavior of $u_1$ and $u_2$ so that when $|\xi_1|$ is large enough, the effect of $u_1$ on $u_2$ is small, and vice versa. Given positive numbers $r_1$ and $a$, if $|\xi_1| \gg 1\,,$ applying the cut and glue-in technique (cf. \S  6) on $u_1$ in a small neighborhood of $\partial B_{\xi_1} (2 r_1)\,,$ and on $u_2$ in a small neighborhood of $\partial B_o (2a)\,,$  one obtains a positive smooth function $u_c$ such that  
$$
u_c = u_1 \ \ \ \ {\mbox{in}} \ \ B_{\xi_1} (r_1) \ \ \ \ {\mbox{and}} \ \ \ \ 
u_c = u_2 \ \ \ \ {\mbox{in}} \ \ B_o (a)\,.
$$
Denote by $K_c$ the $K$-function of $u_c\,.$ [$u_c (x) = u (x)$ and $K_c (x) = K (x)$ for $x$ in a large part of ${\R}^n$, especially when $|x| \gg 1\,.$] The Kelvin transforms of $u_c$ and $K_c$ extend smoothly to  the origin. Furthermore, even for small $\lambda_1\,,$ large $r_1$, and large $a/\lambda_2$, if $\xi_1$ is far away,  we still have 
$$
\sup_{{\R}^n} |K_c - 1| \le 1 - 2^{4/(2-n)} + \epsilon\,,
$$
where the positive constant $\epsilon$ can be chosen to be as small as possible when $|\xi_1| \to \infty\,.$ 
This illustrates the importance of the factor $|\xi_1|$ in (5.23).  On the other hand, if $|\xi_1|/r_1$ is relatively small comparing to $\lambda_2/\lambda_1\,,\,$ then by theorem B, \, $\displaystyle{\sup_{{\R}^n} |K_c - 1|}$ has to be large. We notice the dramatic alternation of $u$ and $u_c$ in the latter case. The sharp change in 
  $K_c$ is mainly introduced by the cut and glue-in process, especially near the regions where $u_1$ or $u_2$ are small [\,cf. (5.26)].


\vspace*{0.3in}

{\bf \Large {\bf 6. \ \ Glue-in the Bubbles}}

\vspace*{0.2in}

In \S 2 the development of bubbles is discussed. We seek to cut and glue a bubble on the solution and keep track of the change in the $K$-function, trying not to disturb it too much. The key point appears to be finding a suitable place to cut and glue. The glue-in argument is based on (2.9) and (2.10). We follow the notations in \S 2.\bk
Let $\delta$ be the constant in (2.9) and $\overline \delta := \delta^{2/(n - 2)}.$ We may assume that $\delta < 1\,.$ Given a positive number $\alpha$ such that $2 ( 1 + \alpha) < n\,,\,$ let $\rho_M$ be a positive number such that 
$$
2 \,{\overline \delta}^{\,{{(n - 2) (n - 2 - 2 \alpha)}\over {2(n + 2)}}} \ge 
\left( {1 \over {1  +  \rho_M^2  }} \right)^{{n - 2}\over 2} \ge \, {\overline \delta\,}^{{(n - 2) (n - 2 - 2 \alpha)}\over {2(n + 2)}} + \,{\overline \delta}^{{n - 2}\over 2}\,. \leqno (6.1)
$$
In particular, 
$$
\left( {1 \over {1  +  \rho_M^2  }} \right) \ge \overline \delta\,.
$$
It follows from (6.1) that  
$$
\left[ \left( {1 \over {1  +  \rho_M^2  }} \right)^{{n - 2}\over 2} - {\overline \delta}^{{n - 2}\over 2}  \right]^{{n + 2}\over {n - 2}} \ge \, {\overline \delta}^{ {{n - 2}\over 2} - \alpha}\,. \leqno (6.2)
$$
For a point $y$ with $|y| = \lambda \,\rho_M$, we have
$$
\left( {\lambda  \over { \lambda^2 +  |y|^2  }} \right)^{{n - 2}\over 2} 
= \left( {\lambda  \over { \lambda^2 +  \lambda^2\,\rho_M^2  }} \right)^{{n - 2}\over 2} 
=  \left( {1 \over {1  +  \rho_M^2  }} \right)^{{n - 2}\over 2} {1\over {\lambda^{{n - 2}\over 2}}}\,. \leqno (6.3)
$$
Here $\lambda$ is the parameter in (2.9) and (2.10).\bk
We may also assume that $R$ is large enough such that $\rho_M < R$ and that $\delta$ is small enough so that, say, $\rho_M > 100$ (in \S 2,  $R$ can be taken arbitrarily large and $\delta$ arbitrarily small). Let $\rho_m := \rho_M - 10$. Take a cut-off function $\phi \in C^\infty ([0, \infty))$ so that $\phi \ge 0$ in $[0, \infty)$, \,$\phi = 1$ in $[0, \, \lambda \, \rho_m)$, \,$\phi = 0$ in $[\lambda \, \rho_M, \,\infty)$, \,$|\phi'| \le C_\phi\, \lambda^{-1}$ and $|\phi'' | \le C_\phi\, \lambda^{-2}$ in $[\lambda\, \rho_m, \, \lambda \, \rho_M]$ for some positive constant $C_\phi$. Let 
$$
w (x) = \phi \,(|x|) \,u_{\lambda\,, o} (x) + [1 - \phi \,(|x|) ] \,u (x_1 + x) \mfor |x| <  \lambda R\,. \leqno (6.4)
$$
We have
\begin{eqnarray*}
\Delta w (x) & = & [ u_{\lambda\,, o} (x) - u (x_1 + x)] \,\Delta \phi (x) + 2 \,[ \btd u_{\lambda\,, o} (x) - \btd_x u (x_1 + x)] \cdot  \btd \phi (x)\\
& \ & \ \ \ \ \, + \, \phi (x) \,\Delta u_{\lambda\,, o} (x) + [1 - \phi (x)]\, \Delta_x u (x_1 + x)\\
& = & [ u_{\lambda\,, o} (x) - u (x_1 + x)] \, \Delta \phi (x) + 2 \,[ \btd u_{\lambda\,, o} (x) - \btd_x u (x_1 + x)] \cdot \btd \phi (x) \\
& \ & \ \ \ \ \, - \, n \,(n - 2) \left [\phi (x) \, u_{\lambda\,, o}^{{n + 2}\over {n - 2}} (x) + [1 - \phi (x)]  \,K (x_1 + x) \,u^{{n + 2}\over {n - 2}} (x_1 + x) \right]
\end{eqnarray*}
for $|x| < \lambda R\,.$ {\it Assume} that 
$$
|K - 1| \le \varepsilon \ \ \ \ {\mbox{in}} \ \ B_{x_1} (\lambda \,R)\,, \leqno (6.5)
$$ 
where $\varepsilon$ is a positive number. 
As $w$ is a non-negative function in $B_o (\lambda \,R)$, there is a function $K_g \in C^\infty (B_o (\lambda \,R))$ (the $K$-function of $w$) such that 
$$
\Delta w + n \,(n - 2) K_g \,w^{{n + 2}\over {n - 2}} = 0 \ \ \ \ {\mbox{in}} \ \ B_o (\lambda \,R)\,. \leqno (6.6)
$$
It follows that
\begin{eqnarray*}
(6.7) \ \  & \  & |K_g (x) - 1|\\
 & = & \bigg\vert {{\Delta w (x)}\over {n (n + 2) w^{{n + 2}\over {n - 2}} (x) }} + 1 \bigg\vert\\
 & \le & 
{{|\Delta \phi (x)|\, |u_{\lambda\,, o} (x) - u (x_1 + x)|}\over { n (n + 2) w^{{n + 2}\over {n - 2}} (x)}} + 
{{2 |\btd \phi (x)|\, |\btd u_{\lambda\,, o} (x) - \btd_x \,u (x_1 + x)|}\over { n (n + 2) w^{{n + 2}\over {n - 2}} (x) }}\\ 
& \ & \ \ \  \,+ \,\Bigg\vert \,{{ u_{\lambda\,, o}^{{n + 2}\over {n - 2}} (x) + [1 - \phi (x)] 
\left[ K (x_1 + x) \,u^{{n + 2}\over {n - 2}} (x_1 + x) - u_{\lambda\,, o}^{{n + 2}\over {n - 2}} (x) \right] }\over
{ \{ u_{\lambda\,, o} (x) + [1 - \phi (x)] 
\left[ u (x_1 + x)  - u_{\lambda\,, o} (x) \right]\, \}^{{n + 2}\over {n - 2}} }} \,- \,1 \,\Bigg\vert   \ \ \ \ \ \ \ \ \ \ \ \ \ \ \ \ \ \ \ \ \   
\end{eqnarray*}
for $x \in B_o (\lambda \,\rho_M) \setminus B_o (\lambda \, \rho_m)\,.$ We have
\begin{eqnarray*}
(6.8) \ \ \ \ \ \ \ \ \ \ \  w^{{n + 2}\over {n - 2}} (x) & = & \left\{ u_{\lambda\,, o} (x) + [1 - \phi (x)] 
\left[ u (x_1 + x)  - u_{\lambda\,, o} (x) \right] \right\}^{{n + 2}\over {n - 2}}\\
&  \ge & \left[ u_{\lambda\,, o} (x) - | u  (x_1 + x) - u_{\lambda\,, o} (x) | \right]^{{n + 2}\over {n - 2}}\\
& \ge & \left[  \left( {1 \over {1  +  \rho_M^2  }} \right)^{{n - 2}\over 2} - 
{\overline \delta}^{{n - 2}\over 2} \right]^{{n + 2}\over {n - 2}} {1\over { \lambda^{{n + 2}\over 2} }}
\, \ge \, {{ {\overline \delta}^{{{n - 2}\over 2}  - \alpha} } \over    { \lambda^{{n + 2}\over 2} }} \ \ \ \ \ \  \ \ \ \ \ \ \ \ \ \ 
\end{eqnarray*}
for $x \in B_o (\lambda \,\rho_M) \setminus B_o (\lambda \, \rho_m)\,.$ As $\phi$ is a radial function,  
$$
\Delta \phi  = {{d^2 \phi}\over {d r^2}} + {{n -1}\over {r}} {{d  \phi}\over {d r }}\,,
$$
where $r = |x|\,.$ Together with the bounds on $\phi'$ and $\phi''$, the fact that $\lambda \rho_m \le |x| < \lambda \rho_M\,,$ and (2.9), (2.10) and (6.8), we obtain 
$$
{{|\Delta \phi (x)| |u_{\lambda\,, o} (x) - u (x_1 + x)|}\over {n \,(n - 2)\, w^{{n + 2}\over {n - 2}} (x) }} \le C \,{\overline \delta}^{\, \alpha} \leqno (6.9)
$$
and
$$
{{2 |\btd \phi (x)| |\btd u_{\lambda\,, o} (x)- \btd_x \, u (x_1 + x) |}\over { n \,(n - 2)\,w^{{n + 2}\over {n - 2}} (x)}} \le C \, {\overline \delta}^{\, \alpha} \leqno (6.10)
$$
for $x \in B_o (\lambda \,\rho_M) \setminus B_o (\lambda \, \rho_m)\,.$ As for the last term in (6.7), we have 

\newpage

\begin{eqnarray*}
(6.11) \  & \ & \Bigg\vert \,{{ u_{\lambda\,, o}^{{n + 2}\over {n - 2}} (x) + [1 - \phi (x)] 
\left[ K (x_1 + x) \,u^{{n + 2}\over {n - 2}} (x_1 + x) - u_{\lambda\,, o}^{{n + 2}\over {n - 2}} (x) \right] }\over
{ \{ u_{\lambda\,, o} (x) + [1 - \phi (x)] 
\left[ u (x_1 + x)  - u_{\lambda\,, o} (x) \right] \}^{{n + 2}\over {n - 2}} }} \,- \,1 \,\Bigg\vert  \\
 & \le & {{ \bigg\vert \, u_{\lambda\,, o}^{{n + 2}\over {n - 2}} (x) - 
\{ u_{\lambda\,, o} (x) + [1 - \phi (x)]  
\left[ u (x_1 + x) - u_{\lambda\,, o} (x) \right] \}^{{n + 2}\over {n - 2}} \bigg\vert }\over {w^{{n + 2}\over {n - 2}} (x) }}\\
& \ & \ \  \,+ \, {{ |K (x + x_1)| \,\bigg\vert \, u^{{n + 2}\over {n - 2}}(x_1 + x) - 
u_{\lambda\,, o}^{{n + 2}\over {n - 2}} (x) \bigg\vert}\over {w^{{n + 2}\over {n - 2}} (x)}} +  |K (x+ x_1) - 1| \,{{u_{\lambda\,, o}^{{n + 2}\over {n - 2}} (x) }\over {w^{{n + 2}\over {n - 2}} (x)}} 
 \ \ \ \  \ \ \ \ \ \ \ \ \ \ \ 
\end{eqnarray*}
for $x \in B_o (\lambda \,\rho_M) \setminus B_o (\lambda \, \rho_m)\,.$  
We note that 
$$
x^{{n + 2}\over {n - 2}} - y^{{n + 2}\over {n - 2}} = {{n + 2}\over {n - 2}} \int^x_y t^{4\over {n - 2}} \, dt \,\le \, {{n + 2}\over {n - 2}} \,x^{4\over {n - 2}}\, |x - y| \mfor x \ge y \ge 0\,. \leqno (6.12)
$$
Also, 
$$
\sup_{B_o (\lambda R)} u_{\lambda\, o}^{4\over {n - 2}} = {1\over {\lambda^2}} \ \ \ \ {\mbox{and}} \ \ \ \ 
\sup_{B_o (\lambda R)} u^{4\over {n - 2}} \le {2\over {\lambda^2}}\,,
$$
as $\lambda$ is small. Hence we obtain
$$
|K (x + x_1)| \,{{ \bigg\vert \, u^{{n + 2}\over {n - 2}} (x_1 + x) - 
u_{\lambda\,, o}^{{n + 2}\over {n - 2}} (x) \bigg\vert}\over {w^{{n + 2}\over {n - 2}} (x)}} 
\le C \lambda^{-2} {{| u (x_1 + x) - u_{\lambda\,, o} (x)|}\over {w^{{n + 2}\over {n - 2}} (x) }} \le C' \, {\overline \delta}^{\, \alpha}  \leqno (6.13)
$$
for $x \in B_o (\lambda \,\rho_M) \setminus B_o (\lambda \, \rho_m)\,.$  Using the left-hand-side inequality of (6.1), we have  
$$
|K (x + x_1)- 1| \,{{u_{\lambda\,, o}^{{n + 2}\over {n - 2}} (x) }\over {w^{{n + 2}\over {n - 2}} (x) }}\le \varepsilon \, {{ \left( {1\over {1 + \rho_m^2}} \right)^{{n + 2}\over 2} \lambda^{-{{n + 2}\over 2}} }\over {{\overline \delta}^{ {{n - 2}\over 2} - \alpha} \lambda^{-{{n + 2}\over 2}}}} \le 2 \,\varepsilon \, {{ \left( {1\over {1 + \rho_M^2}} \right)^{{n + 2}\over 2}   }\over {{\overline \delta}^{ {{n - 2}\over 2} - \alpha}  }} \le 2^{{2n}\over{n - 2}}\, \varepsilon \leqno (6.14)
$$ 
for $x \in B_o (\lambda \,\rho_M) \setminus B_o (\lambda \, \rho_m)\,.$ Here we make use of the fact that $\rho_M \ge 100$ and $\rho_m = \rho_M - 10\,.$ Likewise, using (6.12) and (2.9), we obtain
\begin{eqnarray*}
(6.15) \ \ \ \ \  & \ & {{ \bigg\vert \, u_{\lambda\,, o}^{{n + 2}\over {n - 2}} (x) - 
\{ u_{\lambda\,, o} (x) + [1 - \phi (x)] 
\left[ u (x_1 + x) - u_{\lambda\,, o} (x) \right] \}^{{n + 2}\over {n - 2}} \bigg\vert }\over {w^{{n + 2}\over {n - 2}} (x) }} \ \ \ \ \ \ \ \ \ \ \ \ \ \\
& \le & 
C \lambda^{-2} {{| u (x_1 + x) - u_{\lambda\,, o} (x)|}\over {w^{{n + 2}\over {n - 2}} (x)}}
 \,\le \,C' \,  {\overline \delta}^{\, \alpha} \end{eqnarray*}
for $x \in B_o (\lambda \,\rho_M) \setminus B_o (\lambda \, \rho_m)\,.$ Finally, we glue in the bubble and see that 
$$
u (x_1 + x) = u_{\lambda,\, o} (x) \mfor |x| \le \lambda (\rho_M - 10)\,. 
$$
It follows from (6.5), (6.7), (6.9), (6.10), (6.11), (6.13), (6.14) and (6.15) that 
$$
|K_g - 1|  \le C \, \max \ \{ \varepsilon\,, \  {\overline \delta}^{\,\alpha} \} \leqno (6.16)
$$
in $B_o (\lambda \, \rho_M)\,.$ Here $C$ is a positive constant that depends on $n$, $1 + \varepsilon$ and $C_\phi$ only. For applications in \S 7, we note that $\rho_M$ can be chosen to be large when $\delta$ is small and $R$ is large.


\vspace*{0.3in}

{\bf \Large {\bf 7. \ \ Estimates at the Singularity}}

\vspace*{0.2in}

Let $v$ be a positive smooth solution of the equation 
$$
\Delta v + n (n - 2) \tilde K v^{{n + 2}\over {n - 2}} = 0 \ \ \ \ {\mbox{in}} \ \ {\R}^n. \leqno (7.1)
$$
Consider the Kelvin transform 
$$
u (x) := {{v (y)}\over {|x|^{n - 2}}}\,, \ \ \ \ {\mbox{and}} \ \ \ \ K (x) := \tilde K (y)\,, \ \ \ \ {\mbox{where}} \ \ y := {x\over {|x|^2}}\,, \ \ x \not= 0\,. \leqno (7.2)
$$
We know that $u$ satisfies the equation 
$$
\Delta u + n (n - 2) K u^{{n + 2}\over {n - 2}} = 0  \ \ \ \ {\mbox{in}} \ \ {\R}^n \setminus \{ 0 \}. \leqno (7.3)
$$
We assume that 
$$
\lim_{x \to 0} K (x) = 1 \ \ \ \ {\mbox{and}} \ \ \ \ |x| \cdot |\btd  K (x)| \le C \mfor x \in B_o (1) \setminus \{ 0 \}\,. \leqno (7.4)
$$
It follows from the fact that $K$ is bounded near infinity and the limit in (7.4) that 
$$
|K (x) - 1| \le \kappa^2 \mfor y \in {\R}^n \setminus \{ 0 \}\,, \leqno (7.5)
$$
where $\kappa$ is a positive constant. 
Condition (7.4) is equivalent to 
$$
\lim_{y \to \infty} \tilde K (y) = 1 \ \ \ \ {\mbox{and}} \ \   |y|\cdot|\btd \tilde K (y)| \le C \mfor y \ge 1\,. \leqno (7.6)
$$
\hspace*{0.5in}Assume that $u$ does {\it not} satisfy the slow decay estimate 
$$
u (x) \le C_o |x|^{-{{n - 2}\over 2}}  
$$
for $x$ close to (but not equal to) the isolated singularity at the origin. As discussed in \S 2, bubbles develop and they are described by (2.9) and (2.10).  We choose a small $\delta$ and large $R$ as in \S 2, there is a bubble with center at $\xi_2$ which satisfies (2.9) and (2.10).  Applying the cut and glue-in process in \S 6, we may replace $u$ by another positive smooth 
function which differs from $u$ only on small neighborhood of $\xi_2$, and because of this we still denote the new function by $u$, such that 
$$
u (x) = \left( {{\lambda_2}\over {\lambda_2^2 + |x - \xi_2|^2}} \right)^{{n - 2}\over 2} \ \ \ \ {\mbox{in}} \ \ B_{\xi_2} (a)\,. \leqno (7.7)
$$
Here 
$$
 a = (\rho_M - 10) \,\lambda_2\,, 
$$
where $\rho_M$ is the constant in (6.1). Furthermore, by the choice of $\delta$ and $R$, we have $a \ge \lambda_2\,.$ (Actually we may have $a \gg \lambda_2\,.$)\bk
We choose another set of numbers $\delta'$ (small) and  $R'$ (large). There exists a bubble with center at $\xi_1\,.$ Similar to the above discussion, we may assume that 
$$
u (x) = \left( {{\lambda_1}\over {\lambda_1^2 + |x - \xi_1|^2}} \right)^{{n - 2}\over 2} \ \ \ \ {\mbox{in}} \ \ B_{\xi_1} (r_1)\,, \leqno (7.8)
$$
where 
$$
 r = (\rho_M' - 10) \,\lambda_1\,,  
$$
and $\rho_M'$ is the constant in (6.1) defined by $\delta'\,.$ We also have $r \ge \lambda_1\,.$ 
By choosing $\xi_1$ close to the origin, we may assume that $B_{\xi_2} (a) \cap B_{\xi_1} (r_1) = \emptyset\,.$ We may also assume that  $r_1 /\lambda_1$ and $a/\lambda_2$ are large.  With $\delta$ to be small and $R$ large, and $\xi_1$ and $\xi_2$ close to the origin,  it follows from (6.16) and (7.4)  that (7.5) remains valid.\bk
Under the Kelvin transform with center at the origin and radius $1$, the origin, which is a regular point for both $v$ and $\tilde K$, is sent to infinity. When we apply the Kelvin transform with center at $\xi_2$ and radius $a$ (cf. lemma 5.4), the result is that the infinity is sent to $\xi_2\,.$ It can be seen that the center $\xi_2$ is a regular point for the Kelvin transforms of $u$ and $K$ with center at $\xi_2$ and radius $a$. An argument toward this observation is similar to the proof of lemma 5.4.\bk 
We seek to apply the conclusion in theorem B. As the origin in equation (7.3) is a singularity for $u$, after the Kelvin transformation with center at $\xi_2$ and radius $a$, there is an isolated singularity at $\displaystyle{\,( \xi_2 + a^2 \xi_2/|\xi_2^2|) \in B_{\xi_2} (a)\,.}$ Thus theorem 3.1 may not apply. Consider the following extension of the representation formula with a point singularity.\\[0.2in]
{\bf Theorem 7.9.} \ \ {\it Let $\Omega$ be a bounded domain in ${\R}^n $ with smooth  boundary $\partial \Omega\,.$ For a fixed point $p \in \Omega\,,$ consider a function $u \in C^2 ({\overline \Omega} \setminus \{ p \})\,.$ Assume that there exist positive constants $\mu\,, \ \nu \in (0, 1)\,,$ $C_1$, $C_2$,   
and a small positive number $\delta$ with $B_p (\delta) \in \Omega$, such that 
$$
|\Delta u (x)| \le {{C_1}\over {|x -p|^{n - 1 + \mu} }}\,,\ \ \ \ |\btd u (x)| \le {{C_2}\over {|x-p|^{n - 1 - \nu} }} \ \ \ \ for \ \ x \in B_p (\delta) \setminus \{ p \}\,. \leqno (7.10)
$$
Then for any $\xi \in \Omega \setminus \{ p \}\,,$ we have}
$$
u (\xi) = \int_{\Omega} H (x, \xi) \, \Delta u (x) \, dx + \int_{\partial \Omega} \left[ u (x) \,{{\partial H (x, \xi)}\over {\partial {\bf n}_x}} - H (x, \xi) \,{{\partial u (x)}\over {\partial {\bf n}_x}} \right] \, dS_x\,. \leqno (7.11)
$$

\vspace*{0.2in}


{\bf Proof.} \ \ Given a point $\xi \in \Omega \setminus \{ p \}\,,$ and a small positive number $\varepsilon$ such that $\overline{B_p (\varepsilon)} \cap \overline{B_\xi (\varepsilon)} = \emptyset\,,$ let $\Omega_\varepsilon = \Omega \setminus \{ \overline{B_\xi (\varepsilon)} \cup \overline{B_p (\varepsilon)} \}\,.$ Applying Green's identity we obtain
\begin{eqnarray*}
(7.12) \ \int_{\Omega_\varepsilon} H (x, \xi) \, \Delta u (x) \, dx & = & \int_{\partial \Omega} \left[ 
H (x, \xi)\, {{\partial u (x)}\over {\partial {\bf n}_x}} - u (x) \,{{\partial H (x, \xi)}\over {\partial {\bf n}_x}} \right] \, dS_x\\
& \ & \ \  + \int_{\partial \overline{B_\xi (\varepsilon)}} \left[ H (x, \xi)\, {{\partial u (x)}\over {\partial {\bf n}_x}}  - u (x) \,{{\partial H (x, \xi)}\over {\partial {\bf n}_x}} \right] \, dS_x\\
& \ & \ \ + \int_{\partial \overline{B_p (\varepsilon)}} \left[ H (x, \xi)\, {{\partial u (x)}\over {\partial {\bf n}_x}}  - u (x) \,{{\partial H (x, \xi)}\over {\partial {\bf n}_x}} \right] \, dS_x\,.
\end{eqnarray*}
Here $\bf n$ is the unit outward normal to $\partial \Omega_\varepsilon\,.$ 
It is known that (\cite{Gilbarg-Trudinger}, loc. cit.)
$$
\int_{\partial \overline{B_\xi (\varepsilon)}} \left[ u (x)\, {{\partial H (x, \xi)}\over {\partial {\bf n}_x}} - H (x, \xi)\, {{\partial u (x)}\over {\partial {\bf n}_x}} \right] \, dS_x \to u (\xi) \ \ \ \ {\mbox{as}} \ \ \varepsilon \to 0^+.
$$ 
Furthermore, as $p \not= \xi$, we have 
\begin{eqnarray*}
(7.13) \ \ \ \ \bigg\vert \int_{B_p (\varepsilon) \setminus \{ p \} } H (x, \xi) \, \Delta u (x) \, dx \bigg\vert  & \le & C \int_{B_p (\varepsilon) \setminus \{ p \} }  | \Delta u (x)| \, dx\\ 
& \le & C_1 \int_0^\varepsilon  r^{-\mu} \, dr = {{C_1}\over {1 - \mu}} \varepsilon^{1 - \mu} \to 0 \ \ \ \ \ \ \ \ \ \ \ \ \ 
\end{eqnarray*}
as $\varepsilon \to 0^+.$  It follows from (7.10) that
$$
|u (x)| \le {{C_2}\over {|x-p|^{n - 2 - \nu} }}  \mfor 0 < |x - p| < \delta\,. \leqno (7.14)
$$
Considering (7.10) and (7.14), we have 
\begin{eqnarray*}
& \ & \bigg\vert \int_{\partial \overline{B_p (\varepsilon)}} \left[H (x, \xi)\, {{\partial u (x)}\over {\partial {\bf n}_x}} -  u (x) \,{{\partial H (x, \xi)}\over {\partial {\bf n}_x}}  \right] \, dS_x \bigg\vert\\
& \le & C_3 \int_{\partial \overline{B_p (\varepsilon)}}  \bigg\vert {{\partial u (x)}\over {\partial {\bf n}_x}} \bigg\vert \, dS_x + C_4 \int_{\partial \overline{B_p (\varepsilon)}}  |u (x)| \, dS_x\\
& \le & C_5 (\varepsilon^{1 + \nu} + \varepsilon^\nu) \to 0
\end{eqnarray*}
as $\varepsilon \to 0^+.$ Thus by letting $\varepsilon \to 0^+$ in (7.12) we obtain (7.11).\qed
The {\it main assumption} to apply theorem 7.9 on $v$ in (7.1) is the following. 
$$
|\btd v (y)| \le C |y|^{-\alpha} \mfor |y| \gg 1\,, \leqno (7.15)
$$
where $\alpha$ is a positive constant. 
We also {\it assume, for the moment}, that 
$$
v (y) \ge c \,|y|^{- (n - 2)/2 } \ \ \ \ \mfor |y| \gg 1\,. \leqno (7.16)
$$
It follows from (7.15) that 
$$
v (y) \le C_1 |y|^{1 - \alpha}   \mfor |y| \gg 1\,.\leqno (7.17)
$$ 
Recall that $u (x) = |x|^{2 - n} \,v (x/|x|^2)$ for $x \in {\R}^n \setminus \{ 0 \}\,.$ 
By a direct calculation we obtain
$$
|\btd u (x)|^2 \le C_2 \left[  |x|^{-2 (n - 1)} \, v^2 (y) +   |x|^{-2n}  \,  
|\btd v (y)|^2 \right] \le {C_3 \over {|x|^{2n - 2 \alpha} }} \mfor 0 < |x| < 1\,,
$$
where $y = x/|x|^2$ for $x \not= 0\,.$ From (7.16) we have
$$
u (x) \ge {C_4\over{ |x|^{{n - 2} \over 2} }} \mfor |x| > 0 \ \ {\mbox{small}}\,. \leqno (7.18)
$$
Thus we obtain
$$
\big\vert \btd u^{-4/(n - 2)} (x) \big\vert = \bigg\vert {{\btd u (x)}\over {u^{{n + 2}\over {n - 2}} }}\bigg\vert 
\le {C_5\over{ |x|^{ {{n - 2}\over 2} - \alpha}}} \mfor |x| > 0 \ \ {\mbox{small}}\,. \leqno (7.19) 
$$
Likewise, we have
$$
{{|\btd u (x)|^2}\over {u^{{2n}\over {n - 2}} }} \le {C_6\over{ |x|^{ n - 2 \alpha}}} \mfor |x| > 0 \ \ {\mbox{small}}\,. \leqno (7.20)
$$
It follows from (3.4), (7.3), (7.4) and (7.20) that 
$$
\big\vert \Delta u^{-4/(n - 2)} (x) \big\vert \le {C_7\over{ |x|^{ n - 2 \alpha}}} \mfor |x| > 0 \ \ {\mbox{small}}\,. \leqno (7.21)
$$
We also observe that 
$$
\bigg\vert \btd u^{ - {2\over {n - 2}} } (x) \bigg\vert^2 \le {{C_8}\over {|x|^{n - 2 \alpha}} }
 \mfor |x| > 0 \ \ {\mbox{small}}\,. 
$$
\hspace*{0.5in}As the Kelvin transform with center at $\xi_2$ and radius $a$ does not change the essential features of the singularity at the origin, from (7.19) and (7.21), we may apply theorem 7.9. \bk
We proceed to add a ``base" function to $v$ so that the lower bound (7.16) is fulfilled. 
Let $v_b$ be given by 
$$
v_b (x) =  (|x|^2 + 1)^{(2 - n)/4}  \mfor x \in {\R}^n. \leqno (7.22)
$$ 
The $K$-function of $v_b$ is given by 
$$
K_b (x) = {1\over 2} \left(1 - {{n + 2}\over {2n}} {{|x|^2}\over {|x|^2 + 1}}
\right) 
\mfor x \in {\R}^n. \leqno (7.23)
$$
(See, for example, \cite{Leung.5}, where the constant $n (n - 2)$ is present because of  scaling.) In particular
$$
{{n - 2}\over {4n}} \le K_b (x) \le {1\over 2} \mfor x \in {\R}^n. \leqno
(7.24)
$$
$v_b$ satisfies the equation 
$$
\Delta \,v_b + n (n - 2) K_b \,v_b^{{n + 2}\over {n - 2}} = 0 \ \ \ \ {\mbox{in}} \ \ {\R}^n. 
\leqno (7.25)
$$
We have 
$$
2^{-(n - 2)/4} |x|^{-(n - 2)/2} \le v_b (x) \le  |x|^{-(n - 2)/2} \mfor |x| \ge 1\,,
\leqno (7.26)
$$
and 
$$
|\btd v_b (x)| \le \left( {{n - 2}\over 2} \right) |x|^{-n/2}  \mfor |x| \gg 1\,. 
\leqno (7.27)
$$
(Or else we may use one of the glued-functions similar to those described in \cite{Leung.8}, which is less explicit but with the $K$-function closer to $1\,.$)\bk
Consider 
$$
v_c = v + v_b \ \ \ \ {\mbox{in}} \ \ {\R}^n.
$$
The $K$-function of $v_c$ is given by 
$$
{\tilde K}_c = - {{\Delta v + \Delta v_b}\over {n (n - 2) (v + v_b)^{{n + 2}\over {n - 2}} }} = {{\tilde K v^{{n + 2}\over {n - 2}} + K_b \,v_b^{{n + 2}\over {n - 2}}}\over {(v + v_b)^{{n + 2}\over {n - 2}} }} \ \ \ \ {\mbox{in}} \ \ {\R}^n.
$$
Using (7.24) and the inequalities  
$$
s^{{n + 2}\over {n - 2}} + t^{{n + 2}\over {n - 2}} \le (s + t)^{{n + 2}\over {n - 2}} \ \ \ \ {\mbox{and}} \ \ \ \ 
(s + t)^{{n + 2}\over {n - 2}} \le 2^{4\over {n - 2}} \left(s^{{n + 2}\over {n - 2}} + t^{{n + 2}\over {n - 2}}\right) \leqno (7.28)
$$
for non-negative numbers $s$ and $t$, we obtain
$$
\min\, \left\{ 1 - \kappa^2, \ {{ \min \,\{ (n - 2)/(4n)\,, \ 1 - \kappa^2 \} }\over {2^{4/(n - 2)} }} \right\} \, \le \, {\tilde K}_c \, \le \, \max \, \{ 1 + \kappa^2\,, \ 1/2\} \leqno (7.29)
$$
in ${\R}^n$. Furthermore,
$$
v_c (x) \ge v_b \ge 2^{-(n - 2)/4} \,|x|^{-(n - 2)/2} \mfor |x| \gg 1\,. \leqno (7.30) 
$$
Under the assumption in (7.15), we have  
$$
|\btd v_c (x)| \le |\btd v_b (x)| + |\btd v (x)| \le C |x|^{-\alpha} \mfor |x| \gg 1\,.\leqno (7.31) 
$$
Here we set $0 < \alpha \le n/2\,.$ 
That is, $v_c$ satisfies the conditions (7.4), (7.15) and (7.16). What remains is to show that the essential features of the bubbles are not altered by adding $v_b$ to $v$.\bk 
Let $u_c$, $u$ and $u_b$ be the Kelvin transforms of $v_c$, $v$ and $v_b$, respectively, and $K_c$ the $K$-function of $u_c$, which is given by 
$$
K_c (x) = {\tilde K}_c (x/|x|^2) \mfor x \in {\R}^n \setminus \{ 0 \}\,. \leqno (7.32)
$$
As $u$ does not have slow decay, and $K$ satisfies the conditions in (7.4), bubbles develops according to (2.9) and (2.10). By (7.26) and (7.27), we have
$$
u_b (x) \le |x|^{-(n - 2)/2} \ \ \ \ {\mbox{and}} \ \ \ \ |\btd u_b (x)| \le \left( {{(n^2 - 4) \sqrt n}\over 2} \right)\, |x|^{-n/2} 
\leqno (7.33)
$$
for $|x| > 0$ close to zero. 
From the blow-up analysis in (2.7), (2.9) and (2.10), the distance between $x_o$ and $x_1$ is at most $c_n \,\lambda$, where $c_n$ is a positive constant that depends on $n$ only (cf. \cite{K-M-P-S}). 
In (2.9) and (2.10), we add the bound 
$$
|y| \le \left( {{ d_\varepsilon (x_o)}\over {2 \lambda}} \right) \lambda = {1\over 2} d_\varepsilon (x_o)\,. \leqno (7.34)
$$
Note that, by (2.7),  
$$
{{ d_\varepsilon (x_o)}\over {2 \lambda}} \ge {{M_\varepsilon^{2/(n - 2)} }\over 2}\,,
$$
which can be made arbitrarily large [cf. (2.3)]. By using (7.33) we obtain 
$$
u_b \,(x_1 + y) \le [d_\varepsilon (x_o)/4]^{-(n - 2)/2} \le 2^{n - 2} \,{{\lambda^{{2-n}\over 2} }\over {M_\varepsilon}}  
 \mfor |y| \le {1\over 2} d_\varepsilon (x_o)\,. \leqno (7.35)
$$
Likewise
$$
|\btd u_b \,(x_1 + y)| \le {{2^{n-1} (n-2)}\over { M_\varepsilon^{n/(n - 2)}}}\, \lambda^{-{n\over 2}}  \mfor |y| \le {1\over 2} d_\varepsilon (x_o)\,. \leqno (7.36)
$$
Given a small positive number $\delta$, consider $\delta/2\,.$ Select a large positive number $R$ so that (6.1) if fulfilled. By the blow-up analysis, there exists a positive constant $C$ such that if $M_\varepsilon \ge C$, then (2.9) and (2.10) hold for $\delta/2$ (and hence also for $\delta$). Choose $M_\varepsilon$ to be large enough so that 
$$
{{M_\varepsilon^{2/n-2} }\over 2} \ge R\,, \ \ {\mbox{which \ \ implies}} \ \ \ \ 
{{d_\varepsilon (x_o)}\over 2} \ge R \lambda\,.
$$
From (7.35) and (7.36), if we add the base function to $u$, (2.9) and (2.10) remain valid for $u_c$ with the constant $\delta$, for $|y| < \lambda R \le d_\varepsilon (x_2)/2\,.$\bk
Hence we may assume that $u$ satisfies the lower bound (7.16). By choosing slightly smaller $a$ and $r_1$ if necessary, we may assume that (7.7) and (7.8) are satisfied. 
Thus we are justified to apply theorem 7.9. The arguments in \S 3 and \S 6 lead to similar condition (1.8) for estimate (1.9) to hold. As $\xi_1$ and $\xi_2$ are close to the origin, we have 
$|\xi_1 - \xi_2| \le 1\,.$ 
It follows from (7.29) and (1.9) that there exists a positive constant $C (n, \kappa)$ that depends on $\kappa$ and $n$ only, such that 
$$
{{\lambda_2^2}\over {\lambda_1^2}}  \, \le   \, {{C(n, \kappa)}\over {r_1^4}}\,.\leqno (7.37)
$$
That is,
$$
{{r_1^4}\over {\lambda_1^2}} \le C(n, \kappa) \lambda_2^{-2}. \leqno (7.38)
$$
{\it In order to obtain (7.38), the conditions  are  (7.4) for $K$, and (7.15) for $v$, together with the assumption that $v$ does {\it not} have the slow decay.} 
Fixing $\lambda_2$ in (7.38), we obtain, roughly speaking, that there are no ``deep rooted" bubbles close to the origin. That is, $r_1^2$ cannot be relatively large comparing to $\lambda_1\,.$



\vspace*{0.5in}

{\large D}EPARTMENT OF {\large M}ATHEMATICS\,, \,
{\large N}ATIONAL {\large U}NIVERSITY OF {\large S}INGAPORE\,,\, 
2 {\large S}CIENCE {\large D}RIVE 2\,, \, {\large S}INGAPORE 117543\,, 
{\large R}EPUBLIC OF {\large S}INGAPORE \ \ \ \ 
{\tt matlmc@math.nus.edu.sg}\\[0.1in]
Fax number: 65-779-5452


\newpage

\begin{thebibliography}{30} 
\bibitem{Ambrosetti_et_al} A. Ambrosetti, J. Garcia Azorero and I. Peral, {\it Elliptic variational problems in $R\sp N$ with critical growth,} . J. Differential Equations {\bf 168}  (2000), 10-32. 
\bibitem{Aubin.1} T. Aubin, {\it Sur le probl\'eme de la courbure scalaire prescrite,} Bull. Sci. Math. {\bf 118} (1994), 465-474. 
\bibitem{Caffarelli-Gidas-Spruck} L. Caffarelli, B. Gidas and J. Spruck, {\it
Asymptotic symmetry and local behavior of semilinear elliptic equations with
critical Sobolev growth,} Comm. Pure Appl. Math. {\bf 42} (1989), 271-297.
\bibitem{Chen-Lin.1} C.-C. Chen and C.-S. Lin, {\it On compactness and
completeness of conformal metrics in ${\bf R}^N$,} Asian J. Math. {\bf
1} (1997), 549-559.
\bibitem{Chen-Lin.2} C.-C. Chen and C.-S. Lin, {\it Estimates of the
conformal scalar curvature equation via the method of moving planes,}
Comm. Pure Appl. Math. {\bf 50} (1997), 971-1019.
\bibitem{Chen-Lin.3} C.-C. Chen and C.-S. Lin, {\it Estimates of the
conformal scalar curvature equation via the method of moving planes. II,}
J. Differential Geom. {\bf 49} (1998), 115-178.
\bibitem{Chen-Lin.4} C.-C. Chen and C.-S. Lin, {\it On the asymptotic symmetry of 
singular solutions of the scalar curvature equations,} Math. Ann. {\bf 313} (1999), 229-245.
\bibitem{ChenWX-Li-Ann} W.-X. Chen and C.-M. Li, {\it  A priori estimates for prescribing scalar curvature equations,} Ann. of Math. {\bf 145} (1997), 547-564. 
\bibitem{Cheung-Leung.1} K.-L. Cheung and M.-C. Leung, {\it Asymptotic behavior of  
positive solutions of the equation $\Delta u + K u^{{n + 2}\over {n - 2}} = 0$  
in ${\R}^n$  and positive scalar curvature,\,} Discrete Contin. Dynam. Systems, Added Volume (Proceedings of the International Conference on Dynamical Systems and Differential Equations, Edited by Joshua Du and Shouchuan Hu) (2001), 109-120.
\bibitem{Gidas-Ni-Nirenberg.1}B. Gidas, W.-M. Ni and L. Nirenberg, {\it Symmetry and related properties via the maximum principle,} Comm. Math. Phys. {\bf 68} (1979), 209-243.
\bibitem{Gidas-Ni-Nirenberg.2}B. Gidas, W.-M. Ni and L. Nirenberg, {\it Symmetry of positive solutions of
nonlinear elliptic equations in ${\R}^n$,} Mathematical Analysis and Applications, Part A, p.
369-402, Adv. in Math. Suppl. Stud., 7a, Academic Press, New York-London, 1981.
\bibitem{Gilbarg-Trudinger} D. Gilbarg and N. Trudinger, {\it Elliptic Partial
Differential Equations of Second Order,} 2nd ed., Springer-Verlag, 
Berlin-Heidelberg-New York, 1997.
\bibitem{Kazdan} J. Kazdan,  {\it Prescribing the Curvature of a Riemannian Manifold,}  
CBMS Regional Conference Series in Mathematics, Vol. 57, American Mathematical Society, Providence, Rhode Island, 1985.  
\bibitem{Kato} S. Kato, {\it The scalar curvature equation on open Riemannian manifolds,} Sugaku  Expositions {\bf 14} (2001), 219-236. 
\bibitem{K-M-P-S} N. Korevaar, R. Mazzeo, F. Pacard and R. Schoen,
{\it Refined asymptotics for constant scalar curvature metrics with
isolated singularities,} Invent. Math. {\bf 135} (1999), 233-272.
\bibitem{Leung.1} M.-C. Leung, {\it Conformal scalar curvature equations on
complete manifolds,} Comm. Partial Differential Equations {\bf 20} (1995),
367-417.
\bibitem{Leung.5} M.-C. Leung, {\it Asymptotic behavior of positive solutions
of the equation $\Delta_g u + K u^p = 0$ in a complete Riemannian manifold and
positive scalar curvature,}  Comm. Partial Differential Equations {\bf 24} (1999), 425-462.
\bibitem{Leung.6} M.-C. Leung, {\it Growth estimates on positive solutions 
of the equation $\Delta u + K u^{{n + 2}\over {n - 2}} = 0$ in ${\R}^n,$\,}   
Canad. Math. Bull., {\bf 44} (2001), 210-222. 
\bibitem{Leung.8} M.-C. Leung, {\it Exotic solutions of the conformal scalar curvature equation in $R^n$}, Annales de l'Institut Henri Poincar\'e - Analyse Non Lin\'eaire, 
{\bf 18} (2001), 297-307.   
\bibitem{Leung.7} M.-C. Leung, {\it Conformal scalar curvature equations in open spaces\,,}  Cubo Matem\'atica Educacional, {\bf 3} (2001), 415-443.
\bibitem{Leung.9} M.-C. Leung, {\it Blow up solutions of nonlinear elliptic equations in ${\R}^n$ with critical exponent,} preprint.
\bibitem{YYLi.4}Y.-Y. Li, {\it Fine analysis of blow up and applications,}   Proceedings of the ICCM 1998, Edited by Y. Yang and S.-T. Yau, AMS/IP Studies in Advanced Mathematics, Vol. 20, p. 411-421.
\bibitem{Lin.1}C.-S. Lin, {\it Estimates of the
conformal scalar curvature equation via the method of moving planes III,} Comm. Pure Appl. Math. {\bf 53} (2000), 611-646.
\bibitem{Loewner-Nirenberg}C. Loewner and L. Nirenberg, {\it Partial differential equations invariant under conformal or projective transformations,} Contributions to Analysis (a collection of papers dedicated to Lipman Bers), p. 245-272, Academic Press, New York, 1974.
\bibitem{Mazzeo-Pacard_2} R. Mazzeo and F. Pacard, {\it Constant scalar curvature metrics with isolated singularities,} Duke Math. J. {\bf 99} (1999), 353-418.
\bibitem{M-P-U} R. Mazzeo, D. Pollack and K. Uhlenbeck, {\it Moduli spaces of singular Yamabe metrics,} J. Amer. Math. Soc. {\bf 9} (1996), 303-344.
\bibitem{Melrose} R. Melrose, {\it The Atiyah-Patodi-Singer Index Theorem}, Research Notes
in Mathematics, Vol. 4,  \,A. K. Peters, Ltd., Wellesley, MA, 1993.
\bibitem{Taliaferro.1} S. Taliaferro, {\it On the growth of superharmonic functions 
near an isolated singularity, I.} J. Differential Equations {\bf 158} (1999), 28-47.
\bibitem{Taliaferro.2} S. Taliaferro, {\it Local behavior and global existence of positive solutions of positive solutions of $a u^\lambda \le - \Delta u \le u^\lambda,$\,} 
Annales de l'Institut Henri Poincar\'e - Analyse Non Lin\'eaire, to appear.
\bibitem{Yanagida_Yotsutani} E. Yanagida and S. Yotsutani, {\it Global structure of positive solutions to equations of Matukuma type,} Arch. Rational Mech. Anal. {\bf 134} (1996), 199-226. 


             
         

\end{thebibliography}
\end{document}